\documentclass{amsart}
\usepackage{graphicx}%,psfrag} %pictures package
\usepackage{amsmath,amsthm,amsfonts,amscd,amsxtra,amsopn,amssymb,setspace,verbatim,subfigure}
\usepackage{eucal}
\usepackage{verbatim} %commenting package amongst other features.
\usepackage{ifthen} %conditional package
\usepackage{array}
\usepackage{caption}
\usepackage{stmaryrd}
\def\myPOS#1{\POS}

\def\H{\mathbb H}

\def\Q{\mathbb Q}
\def\R{\mathbb R}

\def\Z{\mathbb Z}
\def\cA{\mathcal A}

\def\cH{\mathcal H}

\def\cL{\mathcal L}

\def\cR{\mathcal R}
\def\cS{\mathcal S}

\def\fg{\mathfrak g}
\newcommand{\into}{\hookrightarrow}

\def\wt{\widetilde}

\def\={:=}
\def\doublearrow#1#2
{\quad\raisebox{.1cm}{$\overset{#1}\lra$}  
\hskip -.67cm
\raisebox{-.1cm}{$\underset {#2}\lra$}\quad}
\def\Map{\operatorname{Map}}

\def\Ker{\operatorname{Ker}}

\def\ind{\operatorname{index}}

\def\nb-{\nobreakdash-}

\def\:{\colon}
\def\iso{\cong}
\def\-1{^{-1}}

\input xy
\xyoption{all}
\numberwithin{equation}{section}

\theoremstyle{plain}
\newtheorem{thm}[equation]{Theorem}
\newtheorem{prop}[equation]{Proposition}
\newtheorem{lemma}[equation]{Lemma}
\newtheorem{cor}[equation]{Corollary}

\theoremstyle{definition}
\newtheorem{defn}[equation]{Definition}
\newtheorem{conj}[equation]{Conjecture}

\newtheorem{question}[equation]{Question}

\theoremstyle{remark}
\newtheorem{rem}[equation]{Remark}

\newcommand\phalf{\frac{p_1}{2}}
\newcommand\dd{\delta}
\newcommand\Dirac{\displaystyle{\not}D}
\newcommand\Aroof{\widehat{A}}

\newcommand\half{\frac{1}{2}}

\newcommand\Witten{\varphi_W}

\newcommand\of{\circ}

\def\Tr{\operatorname{Tr}}

\begin{document}
\title{String Structures and Canonical 3-Forms}
\author{Corbett Redden}
\address{Department of Mathematics, Michigan State University, East Lansing, MI 48824}
\email{redden@math.msu.edu}
\thanks{This work was partially supported by the NSF RTG Grant DMS-0739208.}

%\keywords{string structure, differential characters, positive Ricci curvature, elliptic cohomology}

\begin{abstract}Using basic homotopy constructions, we show that isomorphism classes of string structures on spin bundles are naturally given by certain degree 3 cohomology classes, which we call string classes, on the total space of the bundle.  Using a Hodge isomorphism, we then show that the harmonic representative of a string class gives rise to a canonical 3-form on the base space, refining the associated differential character.  We explicitly calculate this 3-form for homogeneous metrics on 3-spheres, and we discuss how the cohomology theory $tmf$ could potentially encode obstructions to positive Ricci curvature metrics.
\end{abstract}

\maketitle

%\tableofcontents

\section{Introduction}

Degree four characteristic classes arise as obstructions in several ways in math and theoretical physics.  This is analogous to the way the Stiefel--Whitney classes $w_1$ and $w_2$ encode obstructions to orientations and spin structures on a manifold $M$.  One usually encounters the degree four classes when considering structures analogous to the spin structure, but on mapping spaces $\Map(\Sigma, M)$ where $\Sigma$ is a 1 or 2-dimensional manifold.  It is common to say that $\phalf(M)=0\in H^4(M;\Z)$ is the obstruction to forming a {\it string structure} on a manifold $M$.

In this paper, we only deal with a homotopy-theoretic version of string structures.  While geometric notions, such as \cite{CP89, ST04, Wal09}, are necessary for applications, we show we can recover some of this geometric information from the toplogical data for free.  When dealing with these degree 4 classes, one usually must also deal with the associated differential characters.  We naturally obtain globally defined forms representing these characters.  We also speculate on the possibility that the string orientation of $tmf$ may encode obstructions to positive Ricci curvature metrics.  This would be analogous to the obstructions for positive scalar curvature metrics encoded in the spin orientation of $KO$.

Throughout this paper, all manifolds will be compact, connected, oriented, smooth, and without boundary.  We now set up the following notational conventions.  Let $G$ be a compact, simply-connected, simple Lie group and $\lambda\in H^4(BG;\Z)$ a universal characteristic class.  Let $P\overset{\pi}\to M$ be a principal $G$-bundle with connection $\Theta$, and let
\[ \lambda(P)\in H^4(M;\Z), \quad \lambda(\Theta)\in \Omega^4(M), \quad \check{\lambda}(\Theta) \in \check{H}^4(M),\]
be the naturally induced characteristic class, Chern--Weil form, and Cheeger--Simons differential character, respectively.  The differential character is closely related to the Chern--Simons form $CS_{\lambda}(\Theta)\in \Omega^3(P)$.  Finally, $g$ will be a Riemannian metric on $M$.  (In Section \ref{sec:StringStructures}, we also consider more general $G$ and $\lambda$.)

To the characteristic class $\lambda \in H^4(BG;\Z)\iso H^3(G;\Z)\iso \pi_3(G)\iso \Z$, one can associate a topological group $\wt{G}_\lambda$ and homomorphism
\[ \wt{G}_\lambda \to G \]
killing off the corresponding element in $\pi_3(G)$ and inducing isomorphisms in all higher homotopy groups \cite{Sto96, ST04, BSCS07, Hen08, SP09}.  When $G=Spin(k)$ and $\lambda=\phalf$, the resulting group is commonly known as $String(k)$.  While the actual groups $\wt{G}_\lambda$ are not easy to describe, their homotopy type is clearly fixed.  Therefore, we base our constructions only on the homotopy type.  While more concrete models of $\wt{G}_\lambda$ lead to more geometric definitions of $\wt{G}_\lambda$-structures, we only consider the problem of lifting the classifying map from $BG$ to $\wt{BG}_\lambda$.  We call a specific choice of lift a trivialization of the cohomology class $\lambda$.

In Section \ref{sec:StringStructures}, we show that, up to homotopy, such trivializations of $\lambda$ are naturally equivalent to cohomology classes $\cS \in H^3(P;\Z)$ that restrict to $\Omega \lambda \in H^3(G;\Z)$ on the fibers.  Here, $\Omega \lambda$ is the class which universally transgresses to $\lambda$.  These classes $\cS$ are referred to as $\lambda$-trivialization classes.  In the case where $G=Spin(k)$ and $\lambda=\phalf$, we see that the homotopy class of a string structure is equivalent to its string class $\cS$.  An important consequence is that one can describe an element in string bordism by a spin manifold $M$ and string class $\cS \in H^3(Spin(TM);\Z)$.

In fact, these statements hold in greater generality, and Section \ref{sec:StringStructures} considers the more general case where $G$ is a topological group and $\lambda \in H^n(BG;H)$.  Homotopy classes of lifts to $\wt{BG}_\lambda$ still induce canonical classes in $H^{n-1}(P;H)$, and there is an equivalence when $\wt{H}^\bullet(BG;H)=0$ for $\bullet<n$.

In Section \ref{sec:harmonicrepresentatives}, we analyze the harmonic representative of a $\lambda$-trivialization class $\cS$ on $P\overset{\pi}\to M$.  The metric on $P$ is naturally induced by a connection $\Theta$ on $P$ and a Riemannian metric $g$ on $M$.  The harmonic 3-forms on $P$, in an adiabatic limit, were previously analyzed in \cite{Redden08p1}.  In Theorem \ref{thm:canonical3form}, we see that the induced Hodge isomorphism $H^3(P;\R) \overset{\iso}\to \Omega^3(P)$ sends the string class $\cS$ to the 3-form $CS_{\lambda}(\Theta) -\pi^* H_{\cS, g, \Theta}$.  We call the form $H_{\cS, g, \Theta} \in \Omega^3(M)$ the canonical 3-form associated to the $\lambda$-trivialization class, metric, and connection.  Proposition \ref{prop:Hproperties} states
\[ d^*H_{\cS,g,\Theta}=0 \in \Omega^2(M), \text{ and  } \check{H}_{\cS,g,\Theta} = \check{\lambda}(\Theta) \in \check{H}^4(M),\]
where $\check{H}_{\cS,g,\Theta}$ is the induced differential character.  Thus, the form $H_{\cS,g,\Theta}$ lifts $\check{\lambda}(\Theta)$ to take values in $\R$ instead of $\R/\Z$.  This lift is independent of the metric on $M$; the metric picks out the forms with smallest norm lifting $\check{\lambda}(\Theta)$.  

We note that any time one encounters $\check{\lambda}(\Theta)$ and $\lambda(P)=0$, the form $H_{\cS,g,\Theta}$ is relevant because it gives a purely local version of $\check{\lambda}(\Theta)$.  This situation arises in theoretical physics under the guise of anomaly cancellation.  It also arises when constructing the loop group extension bundle $\widehat{LP}\to LP$ restricting to $\widehat{LG}\to LG$.

In Sections \ref{sec:tmf}-\ref{sec:S3}, we deal exclusively with string structures on the frame bundle $Spin(TM)\to M$ of a manifold with spin structure.  Given a string class $\cS$ and Riemannian metric $g$, we use the Levi-Civita connection to produce the the canonical 3-form $H_{\cS,g}$.

Section \ref{sec:tmf} is largely motivational and provides background information on how string structures arise and why they are important.  In particular, we discuss the string orientation
\[ MString \overset{\sigma}\to tmf\]
of the cohomology theory Topological Modular Forms \cite{Hop02} and its tentative relationship to index theory on loop spaces.  This is analogous to the well-understood relationship between $KO$-theory and index theory.  A theorem of Hitchin shows that the spin orientation of $KO$ encodes obstructions to positive scalar curvature metrics.  In the hope of an analogous theorem, Question \ref{question1} asks the following:  If $(M,\cS,g)$ is a closed Riemannian $n$-manifold with string class $\cS$ satisfying both $Ric(g)>0$ and $H_{\cS,g}=0$, does this imply that $\sigma[M,\cS]=0\in tmf^{-n}(pt)?$  

In Section \ref{sec:torsionconnections}, we give an equivalent reformulation of Question \ref{question1}.  One can use the canonical 3-form $H_{\cS,g}$ to modify the Levi-Civita connection, inducing a metric connection $\nabla^{\cS,g}$ with torsion.  Because the metric is used to ``raise an index" of $H_{\cS,g}$, the global rescaling of $M$ determines a canonical 1-parameter family of connections associated to $H_{\cS,g}$.  This converges to the Levi-Civita connection in the large-volume limit.  Proposition \ref{prop:equivalenceofhyp} states that the simultaneous condition $(Ric(g)>0, H_{\cS,g}=0)$ is equivalent to the modified connection having positive Ricci curvature in a small-volume limit.  An interesting sidenote is the alternate description of the Levi-Civita connection.  For a fixed metric $g$, the Levi-Civita connection is the unique metric connection maximizing the Ricci curvature.

In Section \ref{sec:S3}, we examine Question \ref{question1} in the case where $M=S^3$ with a homogeneous metric (under the left or right action of $S^3\iso SU(2)$).  In this case, the answer to Question \ref{question1} is yes, but not if either of the conditions $(Ric(g)>0, H_{\cS,g}=0)$ are weakened in an obvious way.  The only $(\cS,g)$ satisfying both conditions is the string class and round metric induced from $D^4$.  In this case, the string bordism class is obviously 0.  However, there is a 1-parameter family of left-invariant metrics $g$ satisfying $Ric(g)\geq 0$ and $H_{\cR, g}=0$, where $\cR$ is induced by the right-invariant framing.  Since $\sigma[S^3,\cR]=\frac{1}{24}\in tmf^{-3}(pt)$, we see that our question would have a negative answer if one weakens the curvature condition.  Also, one can find Ricci positive metrics $g$ such that $H_{\cR, g}$ is arbitarily small, so one cannot easily weaken the condition $H_{\cS,g}=0$ either.

We close by noting that Section \ref{sec:harmonicrepresentatives} is part of a more general story.  The results of \cite{Redden08p1} and Section \ref{sec:StringStructures} imply that the adiabatic-harmonic representative of a spin$^c$ class gives a canonical 2-form refining the flat differential character $\check{W}_3(\Theta)$.  Similarly, the harmonic representative of an $SU$-class on a $U(n)$-bundle canonically gives a 1-form refining the character $\check{c_1}(\Theta)$.  It appears there is a very general relationship between certain cohomology classes on a bundle $P$, their harmonic representatives, and the associated differential characters.  The author is currently attempting to prove and properly understand these relations.

Much of this paper is based on work done in the author's Ph.D. thesis \cite{Redden06}, and the paper itself was prepared during the author's time at both SUNY Stony Brook and Michigan State University.  The author would like to especially thank Stephan Stolz for encouragement, insight, and numerous suggestions.  He would also like to thank Peter Teichner, Blaine Lawson, Liviu Nicolaescu, Hisham Sati, and Florin Dumitrescu for useful discussions, questions, and comments.

%%%%%%%%%%%%%%%%%%%%%%%%%%%%%%%%%%%%%%%%%%%%%%%%%%%%
%%%%%%%%%%%%%%%%%%%%%%%%%%%%%%%%%%%%%%%%%%%%%%%%%%%%%
\section{Trivializations of characteristic classes}\label{sec:StringStructures}
In this section we make some observations on the general theory of trivializing a characteristic class, and we apply it to the Pontrjagin class $\phalf \in H^4(BSpin;\Z)$ to obtain results in subsequent sections.  In the case of spin structures, or trivializations of $w_2$, the results in this section are quite standard.  In fact, this section is essentially a rewriting of Chapter 2.1 from \cite{LM89} so that it applies in greater generality.

Since we will frequently use the notions of homotopy fibers and Eilenberg--MacLane spaces, we recall a couple of key facts.  If $H$ is an abelian group,\footnote{$H$ is unrelated to the canonical forms in subsequent sections.} then an Eilenberg--MacLane space of type $K(H,n)$ is a space with the only non-trivial homotopy group being $\pi_nK(H,n) \iso H$.  The space $K(H,n)$ is unique up to homotopy and is the classifying space for ordinary cohomology; i.e. for a CW-complex $X$,
\[ H^n(X;H) \iso [X, K(H,n)],\]
where the right-hand side is homotopy classes of based maps $X\to K(H,n)$.  Furthermore, the loopspace functor $\Omega$ induces a homotopy equivalence
\[ \Omega K(H,n) \simeq K(H,n-1)\]
for $n>1$.

The homotopy fiber of a map is defined as the pullback of the path-space fibration.  Given a space $Y$ with basepoint $y_1$, we obtain the pathspace $PY=\{\gamma:[0,1] \to Y | \gamma(1)=y_1\}$.  The natural map $PY \to Y$ given by $\gamma(0)$ is a fibration whose fiber is homotopic to $\Omega Y$.  In fact, $\Omega Y$ acts on the total space of this fibration.  The homotopy fiber $\wt{X}_f$ of a map $X\overset{f}\to Y$ is then the actual pullback of $PY$.  If the homotopy fiber construction is repeated, one obtains a sequence of fibrations homotopic to
\[ \cdots \Omega \wt{X}_f \to \Omega X\overset{\Omega f}\to \Omega Y \to \wt{X}_f \to X\overset{f}\to Y.\]

Now, let $G$ be a connected topological group (of $CW$ type so that standard classifying space constructions apply).  Then, $BG$ is the classifying space for $G$-bundles, and $H^*(BG;H)$ is the cohomology of $BG$ with coefficients in $H$.  The examples we will be concerned with are when $G$ is a classical Lie group such as $SO(n)$ or $Spin(n)$, and $H=\Z$ or $\Z/2.$  Consider a universal characteristic class $\lambda \in H^n(BG;H),$ equivalent to a homotopy class of maps
\[ BG \overset{\lambda}\longrightarrow K(H,n).\]
We fix a specific map $\lambda$ and will not distinguish notationally between the map and the cohomology class.

Let $\widetilde{BG_\lambda}$ be the homotopy fiber of $BG\overset{\lambda}\to K(H,n)$.  This gives rise to the sequence of fibrations up to homotopy
\[ \cdots \to G \overset{\Omega \lambda}\to  K(H,n-1) \to \widetilde{BG_\lambda}\to BG \overset{\lambda}\to K(H,n). \]
Let $P \overset{\pi}\to M$ be a principal $G$-bundle over the space $M$; i.e. $P$ has a free continuous (right) $G$-action with quotient map $\pi : P \to P/G  \iso M.$  Any such bundle $P$ can be obtained as the pullback of the universal bundle 
\[ \xymatrix{ P\ar^{\pi}[d] \ar^{f^*}[r] & EG\ar^{\pi}[d] \\ M \ar^{f}[r] & BG.} \]
Consquently, any $G$-bundle has a natural characteristic class
\[ \lambda(P)\= f^*\lambda \in H^n(M;H).\]

\begin{defn}A trivialization of the characteristic class $\lambda$ on $P$ is a lift of the classifying map to $\widetilde{BG_\lambda}$, i.e. a lift $\widetilde{f}$ $\xymatrix{&\widetilde{BG_\lambda} \ar[d] \\ M \ar@{-->}[ur]^{\widetilde{f}} \ar[r]^{f} &BG}$.\\
We say two trivializations $\wt{f_0}, \wt{f_1}$ are homotopic if they are homotopic through the space of lifts, i.e. if there exists a homotopy $\wt{F}:[0,1]\times M \to \wt{BG_\lambda}$ such that $\wt{F}_{|0} = \wt{f_0}$, $\wt{F}_{|1}=\wt{f_1}$, and  $\wt{F}_{|t}$ is a lift of $f$ for all $t \in [0,1]$.
\end{defn}

\begin{prop}\label{prop:Structures}Let $P \overset{\pi}\to M$ be a $G$-bundle classified by the map $f: M \to BG$. \begin{enumerate}
\item There exists a trivialization of $\lambda$ on $P$ if and only if $\lambda(P)=0\in H^n(M;H)$.
\item If $\lambda(P)=0$, the set of trivializations of $\lambda$ up to homotopy has a free and transitive action of $H^{n-1}(M;H)$; i.e. it is an $H^{n-1}(M;H)$-torsor.
\end{enumerate}
\begin{proof}
Part (1) follows from the definition of the homotopy fiber.  A lift $\wt{f}$ is precisely the choice of a nullhomotopy of $\lambda \circ f :M \to K(H,n)$, and $\lambda \circ f$ is nullhomotopic precisely when the cohomology class $\lambda(P) = 0$.

For part (2), assume an initial trivialization $\wt{f_0}$.  This is equivalent to a global section $\wt{f_0}:M\to f^* \wt{BG_\lambda}$, and $\wt{BG_\lambda} \to BG$ is a fibration with fibers of type $\Omega K(H,n) \simeq K(H,n-1)$.  In fact the H-space $\Omega K(H,n)$ acts fiberwise on $\wt{BG_\lambda}$, so a global section $\wt{f_0}$ induces a fiber homotopy equivalence
\[ \xymatrix{ M \times \Omega K(H,n) \ar[rr]^{\simeq} \ar[dr] & &f^* \wt{BG_\lambda}\ar[dl] \\
& M
}\]
Therefore, the homotopy class of any other section $\wt{f_1}:M\to f^*\wt{BG_\lambda}$ is equivalent to the homotopy class of a function $M\to \Omega K(H,n) \simeq K(H,n-1)$.
\end{proof}\end{prop}

Note that the connectedness of $G$ implies that $BG$ is simply-connected, so we don't have to use local coefficients when dealing with the cohomology of fibers.  The cohomology of any fiber is canonically isomorphic to $H^*(G;H)$, and we have a well-defined ``restriction to fibers" map in cohomology
\[ i^*: H^*(P;H) \to H^*(G;H). \]

\begin{prop}\label{prop:StructureCohClass} \
\begin{enumerate}
\item A trivialization $\wt{f}$ of $\lambda(P)$ gives a canonical cohomology class in $H^{n-1}(P;H)$ that restricts on fibers to the class $\Omega \lambda \in H^{n-1}(G;H)$.
\item The cohomology class in (1) only depends on the homotopy class of $\wt{f}$.
\item Furthermore, $H^{n-1}(M;H)$ acts equivariantly on the homotopy classes of $\lambda$-trivializations and $H^{n-1}(P;H)$ via $\pi^*$.
\end{enumerate}

\begin{proof}
For part (1), consider the universal pullback bundle
\[ \xymatrix{ &EG\ar[d] &\Pi^*EG\ar[l]_{\Pi^*} \ar[d] \\
 &BG &\wt{BG_\lambda}\ar[l]_{\Pi}
}\]
Then, a lift $\wt{f}:M\to \wt{BG_\lambda}$ such that $\Pi \of \wt{f} = f$ is equivalent to a a $G$-equivariant map $\wt{f}^*:P\to \Pi^*EG$ such that $\Pi^* \of \wt{f}^* = f^*$.  

Since $EG$ is contractible, $\Pi^* EG$ is a $K(H,n-1)$ space, as evidenced by the natural homotopy equivalence of fibrations
\[ \xymatrix@R=.25cm{ &G\ar[r] \ar[dd]^{\simeq} & \Pi^*EG \ar[dr] \ar[dd]^{\simeq}\\
&&& \wt{BG}_\lambda. \\
&\Omega BG \ar[r]^{\Omega \lambda} & \Omega K(H,n)  \ar[ur]
}\]
Therefore, any lift $\wt{f}$ is equivalent to $\wt{f^*}:P\to \Pi^*EG \simeq K(H,n-1)$.  When restricted to a fiber, $\wt{f^*}:G\to \Pi^* EG$ is equivalent to $\Omega \lambda: G \to K(H,n-1)$.  This is shown in the following diagram:
\[ \xymatrix{ & & &K(H,n-1)\ar[dl]^{\Pi^*}\ar[d]\\
&P \ar[d] \ar[r]^{f^*} \ar@{-->}@/^/[urr]^>>{\widetilde{f}^*} &EG \ar[d] &\wt{BG_\lambda}\ar[dl]^{\Pi}\\
&M \ar[r]^{f} \ar@{-->}@/^/[urr]^>>{\widetilde{f}} &BG } \]

For part (2), a homotopy $\wt{F}$ between any two trivializations $\wt{f_0}$ and $\wt{f_1}$ naturally lifts to an equivariant homotopy $\wt{F}^*$ between the bundle maps $\wt{f_0}^*$ and $\wt{f_1}^*$.  Therefore, the cohomology class $\wt{f}^* \in H^{n-1}(P;H)$ of a trivialization only depends on the homotopy class of $\wt{f}$.

For part (3), the fiberwise action of $\Omega K(H,n)$ on $\wt{BG_\lambda}\overset{\Pi}\to BG$ naturally pulls back via $\pi^*$ to an action on $\Pi^* EG$.  If $\wt{f_1} = \phi \cdot \wt{f_0}$, where $\phi:M\to \Omega K(H,n)$, then
\[ \wt{f_1}^* = \pi^*\phi \cdot \wt{f_0}^*.\]
Therefore, if two homotopy classes trivializations $[\wt{f_0}], [\wt{f_1}]$ differ by $[\phi] \in H^{n-1}(M,H)$, their natural cohomology classes $[\wt{f_0}^*], [\wt{f_1}^*] \in H^{n-1}(P,H)$ differ by $\pi^*[\phi]$.
\end{proof}\end{prop}

The previous proposition gives a map
\begin{equation}\label{eq:lambdaclass} \{ \lambda \text{-trivializations}\}/\sim \quad \longrightarrow \quad \{ \cS \in H^{n-1}(P;H) \> | \> i^* \cS = \Omega \lambda \in H^{n-1}(G;H) \} \end{equation}
which is equivariant under the natural $H^{n-1}(M;H)$ action.  Here, $\sim$ denotes equivalence up to homotopy.  In general, this map is neither injective nor surjective.  We will refer to such a cohomology class $\cS$ as $\lambda$-trivialization class.

\begin{prop}\label{prop:TrivEquivToCoh}Suppose $\wt{H}^\bullet(G;H)=0$ for $\bullet <n-1$, then map \eqref{eq:lambdaclass} is a bijection.
\begin{proof}
The connectedness of $G$ implies the $E_2$ term in the Leray--Serre cohomology spectral sequence for $EG\to BG$ is 
\[ E^{r,s}_2 \iso H^r(BG;H^s(G;H)),\]
and the contractibility of $EG$ implies that $E^{r,s}_\infty =0$ for $(r,s)\neq (0,0)$.  This, combined with the vanishing of $H^i(G;H)$ for $i < n-1$, implies that the transgression
\[ d_n: E_n^{0,n-1} \iso H^{n-1}(G;H) \longrightarrow E^{n,0}_n \iso H^n(BG;H) \]
is an isomorphism.  In fact, Lemma \ref{lemma:transgression} says that
\[ d_n( \Omega \lambda ) = \lambda.\]

The Leray--Serre cohomology spectral sequence for $P\overset{\pi}\to M$ is pulled back from the sequence for the universal bundle.  This results in the exact sequence
\begin{align*} 0 \to H^{n-1}(M;H) \overset{\pi^*}\to H^{n-1}(P;H) \overset{i^*}\to H^{n-1}(G;H) &\overset{d_n}\to H^n(M;H) \\
\Omega \lambda &\mapsto \lambda(P) \end{align*}
If $\lambda(P)=0$, then the action of $H^{n-1}(M;H)$ is free and transitive on classes in $H^{n-1}(P;H)$ restricting to $\Omega \lambda$.  Since \eqref{eq:lambdaclass} is an equivariant map, and both sides are torsors for $H^{n-1}(M;H)$, it must be a bijection.
\end{proof}
\end{prop}

\begin{lemma}\label{lemma:transgression}Suppose that $\wt{H}^\bullet(X)=0$ for $\bullet<n$.  Then, the cohomology transgression for the pathspace fibration $\Omega X \into PX \to X$ is the inverse of the loop functor; i.e. $d_n^{-1}=\Omega$ in
\[ \xymatrix{ H^{n-1}(\Omega X;H) \ar@/^/[r]^{d_n}_{\iso} & H^n(X;H) \ar@/^/[l]^{\Omega}}. \]
\begin{proof}
For the fibration $\Omega X \into PX \to X$, the transgression and loop functor are related by
\[ \xymatrix{ H^n(X;H) \ar[d] \ar[r]^{\Omega}& H^{n-1}(\Omega X;H) \\
E_n^{n,0} & E_n^{0,n-1} \ar[l]_{d_n} \ar[u]
} \]
This follows from the general relationship between the transgression and cohomology loop suspension \cite{Ser51}.  If $\wt{H}^i(X;H)=0$ for $i<n$, then there is no room for any nontrivial differentials in the Serre spectral sequence until $d_n$.  Therefore, $E_n^{n,0}\iso H^n(X;H)$, $E_n^{0,n-1} \iso H^{n-1}(\Omega X;H)$ and $d_n$ is an isomorphism with inverse $\Omega$.
\end{proof}\end{lemma}

Finally, we wish to make a general note about \textit{stable} cohomology classes.  The usual examples are Chern classes, Pontryagin classes, and Stiefel-Whitney classes, and they correspond to the stable cohomology of classifying spaces for the groups $U(k)$ and $O(k)$.  In general, assume one has a sequence of groups $\{G(k)\}$ and natural inclusions $G(k) \into G(k+1)$ inducing maps
\[ \cdots \to BG(k) \to BG(k+1) \to BG(k+2) \to \cdots \]
such that the cohomology stabilizes.  We then refer to the cohomology of $\displaystyle BG = \lim_{k\to \infty} BG(k)$.  Any cohomology class $\lambda \in H^n(BG;H)$ is stable and defines a sequence of cohomology classes $\lambda_k \in H^n(BG(k);H)$ for all $k$:
\begin{align*}
H^n(BG;H) &\to H^n(BG(k);H)\\
\lambda &\mapsto \lambda_k,
\end{align*}
though the $k$-subscript is usually unnecessary and dropped.  Given a $G(k)$-bundle $P(k)$ classified by $f:M\to BG(k)$, one can stably extend to a $G(k+l)$-bundle $P(k+l)$ by $M\overset{f}\to BG(k) \to BG(k+l)$.  It is obvious that the characteristic class is stable in that $\lambda_{k+l} (P(k+l)) = \lambda_k (P(k)) \in H^n(M;H)$.

\begin{prop}\label{prop:StableTriv}Consider $\lambda \in H^n(BG;H)$.  A trivialization of $\lambda_k$ on any $G(k)$-bundle naturally induces a trivialization of $\lambda$ on any stable extension of $P$.
\begin{proof}
This follows from the naturality of homotopy fibers.  If we denote the inclusion map $\iota:BG(k)\to BG(k+1)$, then 
\[ \wt{BG(k)_\lambda} = \lambda_k^* PK(H,n) = (\lambda_{k+1} \of \iota)^*PK(H,n) = \iota^* \wt {BG(k+1)_\lambda}.\]  Drawing this bundle map, we have
\[ \xymatrix{ &\wt{BG(k)_\lambda}\ar[d]\ar[r] &\wt{BG(k+1)_\lambda} \ar[d]\\
M \ar[ur]^{\wt{f}} \ar[r]^{f} &BG(k) \ar[r]& BG(k+1).
} \]
Any trivialization of $\wt{BG(k)_\lambda}$ naturally extends to a trivialization of $\wt{BG(k+1)_\lambda}$ by composition, and this process can be continued indefinitely.
\end{proof}\end{prop}

To the $G(k_1)$-bundle $P_1\overset{\pi_1}\to M$ and $G(k_2)$-bundle $P_2\overset{\pi_2}\to M$ we can associate the $G(k_1)\times G(k_2)$-bundle $P_1 \times_M P_2 \to M$.  Assuming there are inclusions
\[ BG(k_1)\times BG(k_2) \overset{\iota_1 \times \iota_2}\longrightarrow BG(k_1 + k_2),\]
the bundle $P_1 \times_M P_2$ is also naturally a $G(k_1 + k_2)$-bundle.

Suppose that $H^*(BG;H)=0$ for $*<n$.  The Kunneth formula then implies the additivity of $\lambda \in H^n(BG;H)$: 
\begin{equation}\label{eq:additivity} \lambda (P_1 \times_M P_2) = \lambda(P_1) + \lambda(P_2) \in H^n(M;H).\end{equation}
The bottom square of the below diagram then commutes up to homotopy, implying the existence of the dotted arrow map.
\[ \xymatrix{ \wt{BG(k_1)}_\lambda \times \wt{BG(k_2)}_\lambda \ar@{-->}[r] \ar[d]& \wt{BG(k_1 + k_2)}_\lambda \ar[d]\\
BG(k_1)\times BG(k_2) \ar[d]^{\lambda \times \lambda} \ar[r] & BG(k_1 + k_2) \ar[d]^\lambda \\
K(H,n) \times K(H,n) \ar[r] & K(H,n)
} \]
Therefore, a trivialization of $\lambda$ on the bundles $P_1$ and $P_2$ induces a trivialization of $\lambda$ on $P_1\times_M P_2$ when viewed as a $G(k_1 + k_2)$-bundle (at least up to homotopy).  This can also be seen explicitly in terms of cohomology classes.

\begin{prop}\label{prop:2/3}For $l=1,2$, let $P_l \overset{\pi_l}\to M$ be a $G(k_l)$-bundle.  Let $P_1 \times_M P_2 \to M$ be the $G(k_1)\times G(k_2)$-bundle and $P\to M$ the induced $G(k_1 + k_2)$-bundle.  Assume $H^\bullet(BG(k);H)=0$ for $\bullet<n \>($here $k= k_1, k_2, k_1+k_2$) and $\lambda \in H^n(BG;H)$.  Then up to homotopy, a $\lambda$-trivialization on any two of $\{P, P_1, P_2 \}$ induces a $\lambda$-trivialization on the third.
\begin{proof}
Equation \eqref{eq:additivity} implies the existence of a $\lambda$-trivialization on the third bundle if the other two admit $\lambda$-trivializations.  Proposition \ref{prop:TrivEquivToCoh} states the choice of a trivialization, up to homotopy, is equivalent to a $\lambda$-trivialization class $\cS_i \in H^{n-1}(P_i;H)$ restricting to $\Omega \lambda$ on the fibers.  We now show that the choice of $\lambda$-trivialization class on any two bundles determines one on the third bundle.

Note that there are natural bundle maps
\[ \xymatrix@=.5cm{ & P_1 \times_M P_2 \ar[dl]^{\pi_1} \ar[dr]_{\pi_2} \ar[rr]^{\iota_1\times \iota_2} && P \\
P_1 & & P_2 } \]  
We seek solutions to the equation
\[ (\iota_1\times \iota_2)^* \cS = \pi^* \cS_1 + \pi^* \cS_2.\]
Just as in Proposition \ref{prop:TrivEquivToCoh}, the following commutative diagram is obtained from the Serre spectral sequences for the bundles $P$ and $P_1 \times_M P_2$:
\[ \xymatrix@C-1.5pc{  0 \ar@{->}[r] & H^{n-1}(M)  \ar@{=}[d] \ar[r] &H^{n-1}(P) \ar[r] \ar[d]^{(\iota_1\times \iota_2)^*} &H^{n-1}(G(k_1 + k_2)) \ar[d]^{(\iota_1\times \iota_2)^*} \ar[r] & H^n(M) \ar@{=}[d]\\
0 \ar[r]& H^{n-1}(M) \ar[r] &H^{n-1}(P_1 \times_M P_2) \ar[r] &H^{n-1}(G(k_1)) \oplus H^{n-1}(G(k_2)) \ar[r] & H^n(M) } \]
The cohomology coefficients are all $H$ but suppressed for spacing purposes.

We know that $(\iota_1\times \iota_2)^*\lambda_{k_1+k_2}=\lambda_{k_1} \oplus \lambda_{k_2}$.  For any three classes $\cS, \cS_1, \cS_2$ in the respective bundles, the exact sequence implies
\[ (\iota_1\times \iota_2)^* \cS - \pi_1^* \cS_1 - \pi_2^* \cS_2 = \pi^* \phi\]
for a unique $\phi\in H^{n-1}(M;H)$.  If we fix two of the classes $\cS, \cS_1, \cS_2$, modifying the third by $\phi$ gives us a solution to our desired equation.
\end{proof}\end{prop}

The previous proposition is useful when dealing with cobordism theories.  In the Pontryagin--Thom construction, the relevant extra structure takes place on the stable normal bundle.  Suppose the $m$-manifold $M$ already has a $G$-structure on the stable normal bundle $\nu(M)$.  A lift of the classifying map to $\wt{BG}_\lambda$ induces maps on the Thom spaces, which in turn give an element in the $\wt{G}_\lambda$-bordism group $M\wt{G}_\lambda^{-m}(pt)$.

However, it is often easier or more desirable to describe structures on the tangent bundle.  For any manifold $M$, $TM\oplus \nu(M)$ is canonically isomorphic to the trivial bundle, so Proposition \ref{prop:2/3} often allows us to construct cobordism classes while only dealing with $TM$, or $G(TM)$.

\begin{cor}\label{cor:orientation}Let $\lambda \in H^n(BG;H)$ be a stable class, and suppose $H^\bullet(BG;H)=0$ for $\bullet<n$.  Then, an $m$-manifold $M$ with $G$-structure and $\lambda$-trivialization class $\cS \in H^{n-1}(G(TM);H)$ canonically determines a $\wt{G}_\lambda$-bordism class $[M,\cS]\in M\wt{G}_\lambda^{-m}(pt).$
\end{cor}

We now apply Propositions \ref{prop:Structures}, \ref{prop:StructureCohClass}, \ref{prop:TrivEquivToCoh}, and Corollary \ref{cor:orientation} to recover standard information on spin and spin$^c$ structures as well as a convenient description of string structures.

\subsection{Spin structures}\label{subsec:spin}
For $k>2$, $\pi_1(SO(k)) \iso \Z/2$ and the nontrivial double cover is known as $Spin(k)$.  The Hurewicz image of the generator of $\pi_1(SO(k))$ is the generator of $H^1(SO(k);\Z/2)$, which transgresses to $w_2 \in H^2(BSO(k);\Z/2)$.  It is then clear that 
\[ \wt{BSO(k)}_{w_2} \simeq BSpin(k).\]
Furthermore, there is a spin orientation of $KO$-theory $\alpha:MSpin \to KO$.  Propositions \ref{prop:Structures}, \ref{prop:StructureCohClass}, \ref{prop:TrivEquivToCoh}, and Corollary \ref{cor:orientation} imply the following.

\begin{prop} Let $P\overset{\pi}\to M$ be a principal $SO(k)$-bundle.
\begin{itemize}
\item $P$ admits a spin structure if and only if $w_2(P)=0\in H^2(M;\Z/2)$.
\item The set of spin structures up to isomorphism is naturally equivalent to the set of spin classes $\cS \in H^1(P;\Z/2)$ which restrict to the nontrivial class in $H^1(SO(k);\Z/2)$.
\item The set of spin structures up to isomorphism is a torsor for $H^1(M;\Z/2)$.
\item An oriented $m$-manifold $M$ with spin class $\cS \in H^1(SO(TM);\Z/2)$ gives rise to the bordism class $[M, \cS]\in MSpin^{-m}(pt)$ and $KO$-theory class $\alpha[M,\cS]\in KO^{-m}(pt).$
\end{itemize}\end{prop}

The above statements also have geometric interpretations given by interpreting $H^1(-;\Z/2)$ in terms of double-covers (see Chapter 2.1 of \cite{LM89}).  Then, a spin structure on $P$ is an equivariant double-cover of $P$ restricting fiberwise to the non-trivial double cover of $SO(k)$.

\subsection{Spin$^c$ structures}
For $k>2$, $H^1(SO(k);\Z)=0$, and $H^2(SO(k);\Z) \iso H_1(SO(k);\Z)\iso \pi_1(SO(k)) \iso \Z/2$.  The group $Spin^c(k) = Spin(k)\times_{\Z/2} S^1$ is a non-trivial $S^1$-bundle over $SO(k)$ and hence classified by the generator of $H^2(SO(k);\Z)$; this generator transgresses to $W_3 \in H^3(BSO(k);\Z) \iso \Z/2$.  Therefore, 
\[ \wt{BSO(k)}_{W_3} \simeq BSpin^c(k).\]
Furthermore, there is a spin$^c$ orientation of $K$-theory $MSpin^c \to K$.  Propositions \ref{prop:Structures}, \ref{prop:StructureCohClass}, \ref{prop:TrivEquivToCoh}, and Corollary \ref{cor:orientation} imply the following.

\begin{prop} Let $P\overset{\pi}\to M$ be a principal $SO(k)$-bundle.
\begin{itemize}
\item $P$ admits a spin$^c$ structure if and only if $W_3(P)=0\in H^3(M;\Z)$.
\item The set of spin$^c$ structures up to homotopy is naturally equivalent to the set of classes $\cS \in H^2(P;\Z)$ which restrict to the nontrivial class in $H^2(SO(k);\Z)$.
\item The set of spin$^c$ structures up to homotopy is a torsor for $H^2(M;\Z)$.
\item An oriented $m$-manifold $M$ with spin$^c$ class $\cS \in H^2(SO(TM);\Z)$ gives rise to the bordism class $[M,\cS]\in MSpin^c\>^{-m}(pt)$ and $K$-theory class $\in K^{-m}(pt)$.
\end{itemize}\end{prop}

Again, the above statements all have direct geometric interpretations based on $K(\Z,2) \simeq BS^1$.  A spin$^c$-structure is an equivariant $S^1$-extension of $P$ restricting to the non-trivial extension on fibers.  One can always tensor an $S^1$-bundle over $P$ with the bullback of an $S^1$-bundle on $M$.

\subsection{String structures}\label{subsec:stringstructures}
Let $G$ be any compact simple simply-connected Lie group.  Then, $\pi_2(G)=0$ and $\pi_3(G)\iso H^3(G;\Z)\iso \Z$.  What happens when you kill $\pi_3(G)$?  The 3-connected cover $G\langle 4\rangle \to G$ cannot be a finite-dimensional Lie group, since any connected non-abelian Lie group has non-trivial $\pi_3$.  However, there do exist topological groups $\wt{G} \to G$ which are 3-connected coverings.  Various constructions can be found in \cite{Sto96, ST04, BSCS07, Hen08, SP09}.  The results of all these imply the following (and usually one only needs $G$ semisimple):

Choose a ``level" $\lambda \in H^4(BG;\Z) \iso  H^3(G;\Z).$  Then, there exists a topological group and continuous homomorphism $G\langle \lambda \rangle \to G$ such that $G \langle \lambda\rangle$ has the homotopy type of the fiber of $G \overset{\Omega \lambda}\to K(\Z,3)$.  Applying the classifying space functor gives
\[ BG\langle\lambda \rangle \to BG \overset{\lambda}\to K(\Z,4).\]
When this construction is applied to $G=Spin(k)$ with $\lambda = \phalf \in H^4(BSpin(k);\Z)$, the resulting topological group is known as $String(k)$.  Trivializations of $\phalf$ are commonly referred to as string structures.  Applying the classifying space functor gives us $BString(k)$, and it is clear that
\[ BString(k) \simeq \wt{BSpin(k)}_\phalf.\]

\begin{rem}While multiple models for $G\langle \lambda \rangle$ exist, there is no ``easy" model like the one Clifford algebras provide for the $Spin$ groups.  One must deal with some combination of higher categories, von Neumann algebras, or gerbes, each of which have particular subtleties.  In this paper, we avoid these subtleties by only considering the homotopy type of $G\langle \lambda \rangle$.  While we lose some information, we can characterize lifts of structure groups purely in terms of ordinary cohomology classes.
\end{rem}

\begin{rem}One should be careful when talking about spin, spin$^c$, and string structures up to homotopy.  In addition to ignoring geometric considerations, these structures are naturally categories and have automorphisms; we only deal with isomorphism classes.  The automorphisms play an important role, especially if one wishes to talk about structures locally or glue together manifolds with structures.  See \cite{ST04, Wal09} for more concrete and categorical models of string structures.
\end{rem}

For $k \geq 3$, $Spin(k)$ is simply-connected and compact; therefore $\wt{H}^\bullet(Spin(k);\Z)=0$ for $\bullet <3$.  Also, there is a generalized cohomology theory $tmf$ that has a string orientation $MString \overset{\sigma}\to tmf$.  This orientation is discussed more in Section \ref{sec:tmf}.  Propositions \ref{prop:Structures}, \ref{prop:StructureCohClass}, \ref{prop:TrivEquivToCoh}, and Corollary \ref{cor:orientation} then imply the following statements, which can obviously be rewritten for arbitrary $\lambda \in H^4(BG;\Z)$ (except for the string orientation).

\begin{defn}\label{defn:StringStructure}\label{defn:StringClass}Let $P\overset{\pi}\to M$ be a principal $Spin(k)$-bundle for $k\geq 3$.
\begin{itemize}
\item A string structure on a principal $Spin(k)$-bundle $P\to M$ is a lift of the classifying map to $BString(k)$, i.e. a lift $\widetilde{f}$ $\xymatrix{&BString(k) \ar[d] \\ M \ar@{-->}[ur]^{\widetilde{f}} \ar[r]^{f} &BSpin(k)}$ 
\item A string class $\cS \in H^3(P;\Z)$ is a cohomology class that restricts fiberwise to the stable generator of $H^3(Spin(k);\Z)$.
\end{itemize}
\end{defn}

\begin{prop}\label{prop:StringStructure}Let $P\overset{\pi}\to M$ be a principal $Spin(k)$-bundle for $k\geq 3$.
\begin{itemize} \item $P$ admits a string structure if and only if $\phalf(P)=0 \in H^4(M;\Z).$
\item Up to homotopy, the choice of a string structure is equivalent to the choice of a string class $\cS \in H^3(P;\Z)$.
\item If $\cS$ is a string class, then so is $\cS + \pi^* \phi$ for $\phi \in H^3(M;\Z)$.  This natural action of $H^3(M;\Z)$ on string classes is free and transitive; i.e. the set of string classes is a torsor for $H^3(M;\Z)$.
\item \label{prop:StringNormal}A spin $m$-manifold $M$ with string class $\cS \in H^3(Spin(TM);\Z)$ determines canonical classes $[M,\cS]\in MString^{-m}(pt)$ and $\sigma[M,\cS]\in tmf^{-m}(pt)$.
\end{itemize}\end{prop}

\begin{rem}\label{rem:lowdim}
The cohomology $H^3(Spin(k);\Z)$ does not stabilize until $k=5$, so we briefly describe the stable generator in dimensions 3 and 4.  Under the low-dimensional isomorphisms with the symplectic groups, the groups $Spin(k)$ for $k=3, 4, 5$ are related through the following diagram:
\[ \xymatrix@R=.5cm{ Spin(3) \ar@{^(->}[r] \ar[d]^{\iso} & Spin(4) \ar[d]^{\iso} \ar@{^(->}[r]& Spin(5) \ar[d]^{\iso} \\
Sp(1) \ar@{^(->}^{Id \times Id}[r] \ar@{^(->}[d]& Sp(1)\times Sp(1) \ar@{^(->}[r] \ar@{^(->}[d]& Sp(2) \ar@{^(->}[d] \\
\H  \ar@{^(->}[r]^{Id \times Id}& \H \oplus \H \ar@{^(->}[r]& Gl(\H,2)
} \]
Note that the $Spin(4)$ decomposition is induced by left and right multiplication of the unit quaternions.  The second inclusion $Spin(4)\into Spin(5)$ is isomorphic to the matrix inclusion $Sp(1)\times Sp(1) \into Sp(2)$ along the diagonal.  Since $H^3(Sp(k);\Z)$ stabilizes at $k=1$, we denote by $1$ a generator of $H^3(Sp(1);\Z)\iso H^3(SU(2);\Z)$.  Then, 
\[ \xymatrix@R=.1cm{H^3(Spin(5);\Z) \ar[r]&H^3(Spin(4);\Z) \ar[r] & H^3(Spin(3);\Z) \\
1 \ar@{|->}[r] & (1,1)\ar@{|->}[r] &2 }\]

Furthermore, remember that we originally define $p_1=-c_2$.  Therefore, we see that $\Omega \phalf \in H^3(Spin(3);\Z)$ is twice a generator, and in fact $\Omega \phalf = -2 \Omega c_2 \in H^3(S^3;\Z)$, where $\Omega c_2$ is the usual generator.  
\end{rem}

%%%%%%%%%%%%%%%%%%%%%%%%%%%%%%%%%%%%%%%%%%%%%%%%%%%%
%%%%%%%%%%%%%%%%%%%%%%%%%%%%%%%%%%%%%%%%%%%%%%%%%%%%
\section{Harmonic representative of a string class}\label{sec:harmonicrepresentatives}

In Section \ref{sec:StringStructures} we showed that, up to homotopy, a string structure on a principal $Spin(k)$-bundle $P\to M$ is equivalent to a string class $\cS \in H^3(P;\Z)$.  In this section we consider the harmonic representative of a string class.  This will depend on the choice of a metric on a closed manifold $M$, a connection on $P$, and it involves taking an adiabatic limit.  We also must pass from $\Z$ coefficients to $\R$ coefficients and lose torsion information.  The harmonic representative of $\cS$ is the Chern--Simons 3-form associated to the connection, minus a 3-form on $M$ representing the differential cohomology class $\check{\phalf}(\Theta)$.  The corresponding result holds for arbitrary $\lambda \in H^4(BG;\Z)$, where $G$ is a compact, simple, simply-connected Lie group.

\subsection{Background: Differential Characters}
The canonical 3-form obtained in Theorem \ref{thm:canonical3form} is best understood using the language of differential characters, originally developed in \cite{CS85}.  Another useful reference is \cite{Fre02}.  Let $C_i(M)$ and $Z_i(M)$ denote the group of smooth $i$-chains and cycles on $M$, respectively.  Let $\Omega^i_\Z(M)$ denote the closed differential $i$-forms with integral periods; i.e. their image in $H^i(M;\R)$ lies in the image of $H^i(M;\Z) \to H^i(M;\R)$.  The group of differential characters $\check{H}^i(M)$ is defined as certain homomorphisms satisfying a transgression property:
\begin{align*} \check{H}^i(M) \= \left\{ \parbox{3.2in}{$\chi:Z_{i-1}(M) \to \R/\Z \; | \; \exists \> \omega \in \Omega^i(M) $ satisfying \\
 $\int_\Sigma c^*\omega = \chi(\partial c) \mod \Z \quad \forall c:\Sigma \to M \in C_i(M)$ } \right\} \end{align*}
The form $\omega$ associated to a character $\chi$ must be unique, and in fact $\omega \in \Omega^i_\Z(M)$.  The character $\sigma$ also determines a cohomology class in $H^i(M;\Z)$ whose image in $H^i(M;\R)$ is the same as $[\omega]$.  These two maps induce the following short exact sequences:
\begin{align} 
\label{eq:exact1} 0 \to \frac{\Omega^{i-1}(M)}{\Omega^{i-1}_\Z(M)} \to &\check{H}^i(M) \to H^i(M;\Z) \to 0 \\
\label{eq:exact2} 0 \to H^{i-1}(M;\R/\Z) \to &\check{H}^i(M) \to \Omega^i_\Z(M) \to 0  \\
\label{eq:exact3} 0 \to \frac{H^{i-1}(M;\R)}{H^{i-1}(M;\Z)} \to &\check{H}^i(M) \to H^i(M;\Z) \times_{H^i(M;\R)} \Omega^i_\Z (M) \to 0. \end{align}
In fact, these exact sequences uniquely characterize the groups $\check{H}^i(M)$ \cite{SS08a}; one can refer to $\check{H}^*(M)$ as the differential cohomology of $M$ without specifying the exact model being used, just as one refers to ordinary cohomology without specifying the model.

The importance of differential cohomology is due to the natural factoring of the Chern--Weil homomorphism through $\check{H}^*(M)$.  Any compact Lie group $G$ and universal class $\lambda \in H^{2\bullet}(BG;\Z)$ determine the following for any $G$-bundle $P\to M$ with connection $\Theta$:
\begin{center}
\begin{tabular}{ c l}
Characteristic class  & $\lambda(P) \in H^{2\bullet}(M;\Z)$, \\
Chern--Weil form & $\lambda(\Theta) \in \Omega^{2\bullet}(M)$,\\
Chern--Simons form  & $CS_\lambda(\Theta) \in \Omega^{2\bullet-1}(P)$,\\
Differential character & $\check{\lambda}(\Theta) \in \check{H}^{2\bullet}(M)$.\\
\end{tabular}\end{center}
The integral class and form associated to $\check{\lambda}(\Theta)$ are $\lambda(P)$ and $\lambda(\Theta)$, respectively.

Suppose that $G$ is compact, semisimple, simply-connected and that $\lambda \in H^4(BG;\Z)$.  Then, as discussed in \cite{Fre95}, the associated Chern--Weil form is
\[ \lambda(\Theta) = \langle \Omega \wedge \Omega \rangle \in \Omega^4(M),\]
where $\Omega$ is the curvature of $\Theta$, and $\langle \cdot, \cdot \rangle$ is a suitably normalized $Ad$-invariant inner product on $\fg$.  In this case, the Chern--Simons form is
\[ CS_\lambda(\Theta) = \langle \Theta \wedge \Omega \rangle - \frac{1}{6} \langle \Theta \wedge [\Theta \wedge \Theta] \rangle \in \Omega^3(P), \]
and
\begin{equation}\label{eq:csfiber} [i^* CS_\lambda(\Theta) ] = \Omega \lambda \in H^3(G;\R).\end{equation}
Suppose that $c:X\to M$ is a 3-cycle.  The assumptions on $G$ imply that $c^*P\to X$ admits a global section $p$.  Then,
\[ \check{\lambda}(\Theta)(c) = \int_X p^* (c^*CS_\lambda(\Theta)) \mod \Z.\]
Hence, the information contained in $\check{\lambda}(\Theta) \in \check{H}^4(M)$ is simply the $\R/\Z$-periods of the Chern--Simons 3-form.  One is forced to only consider the $\R/\Z$-periods due to fact that different global sections will give different $\R$-periods.  Note that when $M$ is a connected oriented 3-manifold, \eqref{eq:exact3} implies
\[ \check{H}^4(M) \iso H^3(M;\R)/H^3(M;\Z) \iso \R/\Z,\]
and the isomorphism is given by evaluating on the fundamental cycle $[M]$.  On a 3-manifold, the element $\check{\lambda}(\Theta) \in \check{H}^4(M)\iso \R/\Z$ is often referred to as the Chern--Simons invariant or number of the connection $\Theta$.  First considered in \cite{CS74}, it motivated the theory of differential characters.

\subsection{Hodge isomorphism on $P$}
As a reminder, a Riemannian metric $g$ on an $n$-manifold $M$ induces the Hodge star $*:\Lambda^i TM \to \Lambda^{n-i}TM$, creating the codifferential 
\[ d^*\= (-1)^{n(i+1)+1} * d * :\Omega^\bullet (M)\to \Omega^{\bullet - 1}(M).\]
The Hodge Laplacian is the operator
\[ \Delta_g = dd^* + d^* d = (d+ d^*)^2 : \Omega^\bullet(M) \to \Omega^\bullet(M).\]
When $M$ is closed (compact with no boundary), classical Hodge theory states that there is a canonical isomorphism
\[ H^\bullet(M;\R) \overset{\Pi_{\Ker \Delta_g} }{ \underset{\simeq} \longrightarrow} \Ker \Delta_g \subset \Omega^\bullet (M). \]
We will later denote $\Ker \Delta_g$ by $\cH^\bullet (M)$, though the forms in $\cH^3(P)$ will only be harmonic in a soon to be described limit.

Let $(M,g)$ be a closed Riemannian manifold, and let $P\overset{\pi}\to M$ be a principal $G$-bundle with connection $\Theta$ ($G$ a compact, simple, simply-connected Lie group).  This naturally gives rise to a one-parameter family of right-invariant Riemannian metrics on $P$:
\[ g_\dd \= \dd^{-2}\pi^* g \oplus g_{G}, \quad \dd >0,\]
where $g_G$ is any bi-invariant metric on $G$.  (The metric $g_G$ exists since $G$ is compact, and it is unique up to a scaling constant due to $G$ being simple.)  Conceptually, $g_\dd$ is given by using the connection to decompose $TP$ into horizontal and vertical spaces; the metrics on $M$ and $G$ determine metrics on the horizontal and vertical components respectively.

For any $\dd >0$, we have the harmonic forms $\Ker \Delta^3_{g_\dd} \subset \Omega^3(P).$  In general this finite-dimensional subspace varies with $\dd$, and we will not be concerned with $\Ker \Delta_{g_\dd}$ for any particular $\dd$.  Instead, we analyze the {\it adiabatic limit}, the limit as $\dd \to 0$.  Note that we had to choose the metric $g_G$.  For this reason, it seems natural to introduce the scaling factor $\dd$ and take a limit, thus removing the dependence on the initial choice of $g_G$.  Indeed, this is supported by concrete calculations, where the adiabatic limit appears to be of most interest.  

\begin{thm}[\cite{MM90, Dai91, For95}]\label{thm:extendto0}The 1-parameter space $\Ker \Delta^\bullet_{g_\dd} \subset \Omega^\bullet(P)$ smoothly extends to $\dd=0$.  Furthermore, there is a spectral sequence computing $\displaystyle \lim_{\dd \to 0} \Ker \Delta_{g_\dd}$ which is isomorphic to the Serre spectral sequence.
\end{thm}

The above theorem holds in greater generality, and the context of each cited paper applies to the principal $G$-bundles with metric that we are considering.  The spectral sequence mentioned is a Hodge-theoretic sequence, the details of which are given in \cite{For95} and also summarized in \cite{Redden08p1}.  The fact that $\Ker \Delta_{g_\dd}$ extends continuously to $\dd=0$ (as a path in Grassmannian space) implies that there is still a Hodge isomorphism
\[ H^\bullet(P;\R) \overset{\Pi_{\Ker \Delta_0} }{ \underset{\simeq} \longrightarrow} \lim_{\dd\to 0}\Ker \Delta_{g_\dd} \subset \Omega^\bullet (P). \]
We now introduce the notation
\begin{align*}
\cH^\bullet(M) &\= \Ker \Delta_g \subset \Omega^\bullet(M),\\
\cH^\bullet(P) &\= \lim_{\dd\to 0} \Ker \Delta_{g_\dd} \subset \Omega^\bullet(P), \\
\cH^\bullet(G) &\= \Ker \Delta_{g_G} \subset \Omega^\bullet(G).
\end{align*}
In \cite{Redden08p1}, the spectral sequence interpretation of $\cH^3(P)$ was used to give the following description of harmonic 3-forms on $P$ in the adiabatic limit.

\begin{thm}[Prop 4.5 and Thm 4.6 of \cite{Redden08p1}]\label{thm:harmonic3forms}Consider $(P\overset{\pi}\to M, g, \Theta)$ where $G$ is a compact simple Lie group.  If $\lambda(P)=0 \in H^4(M;\R)$, then
\[ \cH^3(P) = \R[ CS_\lambda(\Theta) - \pi^* h] \> \oplus \> \pi^* \cH^3(M),\]
where $h \in \Omega^3(M)$ is the unique coexact form satisfying $dh = \lambda(\Theta)$.
\end{thm}

When $G$ is also simply-connected, the Serre spectral sequence gives the following exact sequence, as seen in Proposition \ref{prop:TrivEquivToCoh}:
\[ \xymatrix@R=.1cm{0 \ar[r] &H^3(M;\Z) \ar[r]^{\pi^*}&H^3(P;\Z) \ar[r]^{i^*} &H^3(G;\Z) \ar[r]^{d_4}&H^4(M;\Z) \\
& & \cS \ar@{|-->}[r]^{?} & \Omega \lambda \ar@{|->}[r] &\lambda(P) }\]

\begin{thm}\label{thm:canonical3form}Consider $(P\overset{\pi}\to M, g, \Theta)$ where $G$ is a simply-connected compact simple Lie group.  Suppose that $\lambda(P)=0\in H^4(M;\Z)$ and that $\cS \in H^3(P;\Z)$ is a $\lambda$-trivialization class, i.e. $i^* \cS = \Omega \lambda \in H^3(G;\Z)$.  Then, the image of $\cS$ under the Hodge isomorphism is of the form
\begin{align*}
H^3(P;\Z) &\to H^3(P;\R) \overset{\Pi_{\Ker \Delta_0}} {\underset{\simeq} \longrightarrow} \cH^3(P) \subset \Omega^3(P)
\\
\cS &\longmapsto \qquad \quad  CS_\lambda(\Theta) - \pi^* H_{\cS, g, \Theta},
\end{align*}
where $H_{\cS, g, \Theta} \in \Omega^3(M)$.  Alternatively, $\Pi_{\Ker \Delta_0} \cS - CS_\lambda(\Theta) \in \pi^* \Omega^3(M)$.
\begin{proof}
We use the description of $\cH^3(P)$ given in Theorem \ref{thm:harmonic3forms}.  The orthogonal decomposition of $\cH^3(P)$ in Theorem \ref{thm:harmonic3forms} corresponds to a splitting $H^3(P;\R) \iso H^3(G;\R) \oplus H^3(M;\R)$.  We know $\pi^*H^3(M;\R)$ restricts to $0 \subset H^3(G;\R)$.  Both $CS_\lambda(\Theta)-\pi^*h$ and $\cS$ cohomologically restrict to $\Omega \lambda \in H^3(G;\R)$, as mentioned in \eqref{eq:csfiber}, so 
\[ \Pi_{\Ker \Delta_0} \cS - (CS_\lambda(\Theta) - \pi^*h) \in \pi^* \cH^3(M).\]
Therefore, the harmonic representative of $\cS$ must be of the form
\[ CS_\lambda(\Theta) - \pi^*h - \pi^*h',\]
with $h' \in \cH^3(M)$, and we define 
\[ H_{\cS, g, \Theta} \= h + h' \in \Omega^3(M).\]
\end{proof}\end{thm}

\begin{rem}The above theorem does not hold without taking an adiabatic limit.  For a general $\dd>0$,
\[ \Pi_{\Ker \Delta_{g_\dd}} \cS - CS_\lambda(\Theta) \notin \pi^* \Omega^3(M),\]
but instead will contain forms with bi-degree (2,1) and (1,2) in the (horizontal, vertical) decomposition of $\Omega^3(P)$ given by the connection.
\end{rem}

\begin{rem}When restricted to the fibers, the Chern--Simons form is the standard harmonic (bi-invariant) form representing $\Omega \lambda$; i.e. $i^* CS_\lambda(\Theta) \in \cH^3(G)$.  Just as $\cS$ is a cohomological extension of $\Omega \lambda$ to all of $P$, we see that $\Pi_{\Ker \Delta_0} \cS = CS_\lambda(\Theta)-\pi^*H_{\cS,g,\Theta}$ is a harmonic extension of $\Omega \lambda$ to all of $P$.
\end{rem}

\subsection{Properties of canonical 3-form}
Theorem \ref{thm:canonical3form} gives a canonical construction
\begin{align}\label{eq:construction}
 \{ \text{$\lambda$-triv classes} \} \times Met(M) \times \cA(P) &\longrightarrow \Omega^3(M) \\
  \cS, g, \Theta &\longmapsto H_{\cS, g, \Theta}. \nonumber \end{align}
We call $H_{\cS, g, \Theta}$ the canonical 3-form associated to $(\cS, g, \Theta)$.  While Theorem \ref{thm:canonical3form} only uses information about $\cS$ as a class in $H^3(P;\R)$, the integrality becomes necessary when understanding $H_{\cS,g,\Theta}$ in terms of differential characters.  The exact sequence \eqref{eq:exact1} gives rise to
\begin{align}\label{eq:exact4} 0 \to \Omega_\Z^3(M) \to \Omega^3(M) &\to \check{H}^4(M) \to H^4(M;\Z) \to 0 \\
H_{\cS, g, \Theta} &\mapsto \check{H}_{\cS, g, \Theta} \nonumber \end{align}
where the character $\check{H}_{\cS, g, \Theta}$ obtained via $\Omega^3(M) \to \check{H}^4(M)$ is given by simply by integrating $H_{\cS, g, \Theta}$ on cycles and reducing mod $\Z$.  Also, note that $H^3(M;\Z)$ acts naturally on $\{ \lambda$-triv classes$\}$, and it also acts on $\Omega^3(M) \times Met(M)$ by adding a harmonic representative.

\begin{prop}\label{prop:Hproperties}The construction \eqref{eq:construction} is equivariant with respect to the natural action of $H^3(M;\Z)$; i.e. $H_{\cS + \pi^* \phi, g, \Theta} = H_{\cS, g, \Theta} + \Pi_{\Ker \Delta_g} \phi.$\\
Furthermore, the forms $H_{\cS,g,\Theta}$ satisfy the following:
\begin{itemize}
\item $d^*H_{\cS,g,\Theta}=0 \in \Omega^2(M)$,
\item $dH_{\cS, g, \Theta} = \lambda(\Theta) \in \Omega^4(M),$
\item $\check{H}_{\cS, g, \Theta} = \check{\lambda}(\Theta) \in \check{H}^4(M).$
\end{itemize}
\begin{proof}
The action of $H^3(M;\Z)$ on $\lambda$-trivialization classes is given by addition under $\pi^*$, and the action on $\Omega^3(M)$ is given by adding the harmonic representative (with respect to a fixed metric $M$).  Theorem \ref{thm:harmonic3forms} implies that for $\phi \in H^3(M;\Z)$ with harmonic representative $\Pi_{\Ker \Delta_g}\phi \in \cH^3(M)$,
\[ \pi^*( \Pi_{\Ker \Delta_g} \phi )= \Pi_{\Ker \Delta_0} (\pi^*\phi) \in \Omega^3(P).\]

The property $d^*H_{\cS,g,\Theta}=0$ also follows directly from Theorem \ref{thm:harmonic3forms}.  That $\Pi_{\Ker \Delta_0}\cS$ is closed implies
\begin{align*} d\left( CS_\lambda(\Theta) - \pi^* H_{\cS,g,\Theta} \right) &= 0 \\
\pi^* \lambda(\Theta) - \pi^* dH_{\cS,g,\Theta} &=0\\
dH_{\cS, g, \Theta} &= \lambda(\Theta),
\end{align*}
with the last equality following from $\pi^*$ being injective on forms.

Finally, suppose that $X \overset{c}\to M$ is a smooth 3-cycle on $M$.  Then, the value of $\check{\lambda}(\Theta)$ on $(X, c)$ is 
\[ \check{\lambda}(\Theta)(c) = c^* \check{\lambda}(\Theta) \in \check{H}^4(X) \iso \R/\Z.\]
Standard obstruction theory implies that $c^*P \to X$ admits a global section $p:X \to c^*P$, and it is easy to see that 
\[ p^* c^* \check{CS_\lambda}(\Theta) = p^* \check{CS_\lambda}(c^*\Theta) = c^* \check{\lambda}(\Theta) \in \check{H}^4(X).\]
Because $CS_\lambda(\Theta) - \pi^*H_{\cS,g,\Theta} \in \Omega^3_\Z(P)$, then
\[ \check{CS_\lambda}(\Theta) = \pi^*\check{H}_{\cS,g,\Theta} \in \check{H}^3(P) \]
and hence
\[ c^*\check{\lambda}(\Theta) = p^* c^* \check{CS_\lambda}(\Theta) = p^* c^* \pi^* \check{H}_{\cS, g,\Theta} = p^* \pi^* c^* \check{ H }_{\cS, g,\Theta}= c^*\check{H}_{\cS, g,\Theta} \in \check{H}^4(X).\]
This implies that for all 3-cycles $X\overset{c}\to M$
\[ \check{H}_{\cS,g,\Theta} (c) = \check{\lambda}(\Theta) (c) ,\]
and hence $\check{\lambda}(\Theta) = \check{H}_{\cS,g,\Theta} \in \check{H}^4(M).$  (This also implies $dH_{\cS,g,\Theta}=\lambda(\Theta)$.)
\end{proof}\end{prop}

Integrating the form $H_{\cS, g,\Theta}$ naturally gives values in $\R$, and Proposition \ref{prop:Hproperties} says that reducing mod $\Z$ gives the same values as $\check{\lambda}(P)$.  In other words, the choice of a $\lambda$-trivialization class naturally gives a lift
\begin{equation}\label{eq:lifta} \xymatrix{ &\R \ar[d]\\
Z_3(M) \ar[r]^{\check{\lambda}(\Theta)} \ar[ur]^{H_{\cS, g,\Theta}}& \R/\Z,
} \end{equation}
and the action of $H^3(M;\Z)$ modifies the lift by the induced map $Z_3(M) \to \Z$.  While the actual form $H_{\cS, g, \Theta}$ depends on the choice of a metric, the above lift does not.

\begin{prop}The lift in \eqref{eq:lifta} is independent of the choice of metric $g$.
\begin{proof}If $g_0, g_1$ are two different metrics, then \eqref{eq:exact4} implies
\[ H_{\cS, g_1, \Theta} - H_{\cS, g_0, \Theta} \in \Omega^3_\Z(M).\]
The space of Riemannian metrics is contractible, so
\[ H_{\cS, g_1, \Theta} - H_{\cS, g_0, \Theta} \in d\Omega^2(M).\]
\end{proof}\end{prop}

The role of the metric in the construction \eqref{eq:construction} is to pick out the forms $H_{\cS, g,\Theta}$ with smallest norm still satisfying $\check{H}_{\cS, g,\Theta} = \check{\lambda}(\Theta)$.  We denote the lift by $H_{\cS,\Theta}\in \check{H}^4_\R(M)$.  Here, we use the non-standard notation of $\check{H}^4_\R(M)$ to denote characters $Z_3(M) \to \R$ satisfying the usual transgression assumption.

To summarize, the construction \eqref{eq:construction} induces lifts of the standard differential character construction, which are encoded in the following diagram:
\[ \xymatrix{ \{ \lambda\text{-triv classes} \} \times Met(M) \times \cA(P) \ar[rr]^<<<<<<{H_{\cS, g,\Theta}}\ar[d]&&\Omega^3(M) \ar[d]\\
\{ \lambda \text{-triv classes} \} \times \cA(P) \ar[rr]^{H_{\cS,\Theta}} \ar[d]&&\check{H}_\R^4(M) \ar[d] \\
\cA(P) \ar[rr]^{\check{\lambda}(\Theta)} && \check{H}^4(M).
} \]

\begin{rem}In \cite{ST04}, Stolz and Teichner define a {\it geometric} trivialization of $\lambda(P)$ as a trivialization of the extended Chern--Simons field theory on $P\to M$.  This includes defining a lift of the differential character $\check{\lambda}(\Theta)$ to take values in $\R$, and it aligns nicely with the above construction above.  In fact, if $H^3(M;\Z)$ has no torsion, then the choice of a lift of $\check{\lambda}(\Theta)$ to $\check{H}^4_\R(M)$ is equivalent to the choice of a $\lambda$-trivialization class.  In \cite{Wal09}, an explicit model for string structures is given in terms of trivializations of a Chern--Simons 2-gerbe, and it is shown that a string structure produces a 3-form on $M$.  The 3-forms obtained in our construction are a proper subset of those obtained in \cite{Wal09}, analogous to the relationship between forms representing a de Rham class and harmonic forms.
\end{rem}

Note that one can also directly define the lift $H_{\cS,\Theta}$ without using the Hodge isomorphism.  On a 3-cycle $c:X\to M$,
\begin{equation}\label{eq:lift} H_{\cS, \Theta}(c) = \int_X p^*\left( CS_\lambda(c^*\Theta) - c^*S \right) \end{equation}
where $p$ is any global section, and $S$ is any de Rham representative of $\cS$.  This is a simple consequence of $S=CS_\lambda(\Theta) - H_{\cS,g,\Theta} + d\beta$.  It is also easy to verify that the integral on the right-hand side is independent of $p$.  In cases like the following Lemma \ref{lem:trivbundle}, this allows us to calculate the form $H_{\cS, g,\Theta}$ without solving a differential equation.

Suppose the $G$-bundle $P\overset{\pi}\to M$ is topologically trivial.  Then, the choice of a global section $p:M\to P$ is equivalent to a trivialization $P \iso M\times G$.  The canonical $\lambda$-trivialization on $M\times G$ induces one on $P$, and the corresponding cohomology class is given by the Kunneth isomorphism
\begin{align}\label{eq:trivbundleclass} H^3(P;\Z) &\iso H^3(M;\Z) \oplus H^3(G;\Z) \\
\cS &\leftrightarrow (0,\Omega \lambda). \nonumber
 \end{align}

\begin{lemma}\label{lem:trivbundle}Suppose $P\overset{\pi}\to M$ is a trivial bundle with $\lambda$-trivialization class $\cS$ induced by the trivialization $p:M\to P$.  Then,
\[ H_{\cS,g,\Theta} - p^* CS_\lambda(\Theta) \in d\Omega^2(M).\]
In particular, $p^*CS_\lambda (\Theta) = H_{\cS,\Theta}$ as elements of $\check{H}^4_\R(M)$.  If $d^* p^* CS_\lambda(\Theta)=0$, then $p^*CS_\lambda(\Theta)=H_{\cS,g,\Theta}\in \Omega^3(M)$.
\begin{proof}
As seen in \eqref{eq:trivbundleclass}, $p^* \cS =0 \in H^3(M;\Z)$.  Therefore, \eqref{eq:lift} simplifies to 
\[ \int_X c^* H_{\cS,g,\Theta} = \int_X c^* p^* CS_\lambda(\Theta)\]
for all 3-cycles $c:X\to M$, so $[H_{\cS,g,\Theta} - p^* CS_\lambda(\Theta)]=0\in H^3(M;\R)$.
\end{proof}\end{lemma}

One usually chooses $\lambda \in H^4(BG;\Z) \iso \Z$ to be the generator.  This is because $\wt{BG}_\lambda$ is the universal extension.  This universality is also reflected in the associated canonical 3-forms.

\begin{prop}If $\cS \in H^3(P;\Z)$ is a $\lambda$-trivialization class and $\ell \in \Z$, then $\ell \cS$ is an $\ell \lambda$-trivialization class, and
\[ H_{\ell\cS, g, \Theta} = \ell H_{\cS, g, \Theta} \in \Omega^3(P).\]
\begin{proof}The first statement is obvious, and the second follows from the linearity of the Hodge isomorphism.
\end{proof}\end{prop}

We now apply the above construction to $G=Spin(k)$ for $k\geq 3$ with $\lambda = \phalf \in H^4(BSpin(k);\Z)$ to canonically produce 3-forms associated to string structures.  Since $Spin(4)\iso SU(2)\times SU(2)$ is not simple, we define the canonical 3-form when $k=4$ to be the one obtained by stabilizing to $Spin(5)$, a process that does not affect $\check{\phalf}(\Theta)$.\footnote{The arguments in \cite{Redden08p1} can be extended to semisimple groups, and Theorem \ref{thm:stringforms} also holds for $k=4$ without stabilizing.}

\begin{thm}\label{thm:stringforms}Let $P\overset{\pi}\to M$ be a principal $Spin(k)$-bundle ($k \geq 3$) with connection $\Theta$ over the Riemannian manifold $(M,g)$.  Under the Hodge isomorphism (in an adiabatic limit), a string class $\cS \in H^3(P;\Z)$ is represented by 
\[ \Pi_{\Ker \Delta_0} \cS = CS_\phalf(\Theta) - \pi^* H_{\cS,g,\Theta} \in \Omega^3(P).\]
Furthermore, the canonical form $H_{\cS,g,\Theta} \in \Omega^3(M)$ satisfies
\begin{itemize} \item $d^*H_{\cS, g, \Theta} = 0 \in \Omega^2(M),$
\item $\check{H}_{\cS, g, \Theta} = \check{\phalf}(\Theta) \in \check{H}^4(M),$
\item the construction of $H_{\cS, g, \Theta}$ is equivariant with respect to the natural action of $H^3(M;\Z)$.
\end{itemize}\end{thm}

In particular, consider the case where $(M,g)$ is a Riemannian manifold with spin structure satisfying $\phalf(M)=0\in H^4(M;\Z)$.  Then, we can let $P=Spin(TM)$, and we call a string structure on $Spin(TM)$ a string structure on $M$.  Letting $\Theta$ be the Levi-Civita connection, this gives a map
\begin{align}\label{eq:Riemannianconstruction}
\{\text{String classes on }M\} \times Met(M) &\longrightarrow \Omega^3(M) \\
\cS, g&\longmapsto H_{\cS, g} \nonumber
\end{align}

%%%%%%%%%%%%%%%%%%%%%%%%%%%%%%%%%%%%%%%%%%%%%%%%%%%%%%
%%%%%%%%%%%%%%%%%%%%%%%%%%%%%%%%%%%%%%%%%%%%%%%%%%%%%%
\section{Canonical 3-forms and the string orientation of $tmf$}\label{sec:tmf}
We now discuss some background on how string structures arise and give a possible new application of the canonical 3-forms $H_{\cS,g}$ from \eqref{eq:Riemannianconstruction}.  To do so, let us first recall some classical results from index theory (an excellent source is \cite{LM89}).  Suppose $M$ is an oriented closed manifold.  A priori, one cannot form a spinor bundle $SO(M)\times_{SO(n)} S^\pm \to M$, because the spinor representations $SO(n)\to Gl(S^\pm)$ are only projective.  The choice of a spin structure, previously discussed in Section \ref{subsec:spin}, allows one to define the spinor bundle $S^{\pm}_M \= Spin(M)\times_{Spin(n)} S^\pm$ and Dirac operator $\Dirac_M:\Gamma(S^\pm)\to \Gamma(S^\mp)$.

While the Fredholm operator $\Dirac_M$ depends on the spin structure, the Atiyah--Singer index theorem states that its index does not, and in fact 
\[ \ind (\Dirac_M) = \Aroof(M) \in \Z. \]
Here, $\Aroof(M)$ is a topological invariant determined by a manifold's Pontryagin classes and is defined for any oriented manifold.  In general $\Aroof(M) \in \Q$, but $\Aroof(M)\in \Z$ when $w_2(M)=0$.  There is also a refinement of $\Aroof(M)$ given by the spin-orientation $\alpha:MSpin\to KO$.  This refinement can be thought of as the Clifford-linear index, and it does depend on the spin structure.
\begin{equation}\label{eq:spinindex} \xymatrix{& KO^{-n}(pt) \ar[d] \\
MSpin^{-n}(pt) \ar[r]^{\Aroof} \ar[ur]^{\alpha} & \Z .} \end{equation}
The $KO$-invariants usually appear in family index theorems, but they also contain interesting information for a single manifold due to the torsion in $KO^{-*}(pt)$.

Index theory is now a central part of mathematics, and one of its powerful applications is to the problem of when a closed manifold admits positive scalar curvature metrics.  The Lichnerowicz--Weitzenb\"ock formula $\displaystyle \Dirac_M^2 = \nabla^*\nabla + \frac{s}{4}$, which relates $\Dirac_M$ to a positive operator and the scalar curvature $s$, implies the following:  If a closed spin manifold $M$ admits a metric of positive scalar curvature, then $\ind (\Dirac_M)=\Aroof(M)=0$ \cite{Lic62}.  Furthermore, $\alpha[M]=0\in KO^{-n}(pt)$ for all spin structures \cite{Hit74}.  In fact, for simply-connected spin manifolds of dimension $\geq 5$, all the $\alpha$-invariants vanish if and only if $M$ admits a metric of positive scalar curvature \cite{Sto92}.

There is an analogous story, though not fully developed, involving the Witten genus, index theory on loop spaces, and elliptic cohomology.  In \cite{Wit88}, Witten used intuition from theoretical physics and defined a topological invariant $\Witten (M)$ known as the Witten genus.  He claimed it should be the $S^1$-equivariant index of the Dirac operator on the free loop space $LM$; i.e.
\[ ``\ind^{S^1} \Dirac_{LM}" = \Witten (M).\]
We place the left-hand side in quotes because of analytic difficulties in defining a good theory of Fredholm operators on infinite-dimensional manifolds.  However, the Witten genus (and other elliptic genera) are well-defined, and one can make formal sense of index theory on $LM$ by using localization formulae or the representation theory of loop groups.  For a good overview on these ideas, see \cite{Liu96}.  While $\Witten (M) \in \Q\llbracket q \rrbracket [q^{-1}]$ for any oriented manifold, for a string manifold $\Witten (M)$ is the $q$-expansion of a modular form ($MF$) with integer coefficients and weight $n/2$, and we say $\Witten (M) \in MF_n$.  The intuitive reason is that when $\phalf(M)=0$, one can define the spinor bundle on $LM$ \cite{CP89}.  We wish to form $LSpin(M) \times_{LSpin(n)}S \to LM$ where $S$ is a positive-energy represenation of $LSpin(n)$.  However, these representations are all projective, so one must pass to an $S^1$-extension $\widehat{LSpin(n)} \to LSpin(n)$.  Topologically, our string class $\cS \in H^3(Spin(M);\Z)$ transgresses to a class in $H^2(LSpin(M);\Z)$ that defines an isomorphism class of $S^1$-extension $\widehat{LSpin(M)}\to LSpin(M) \to LM$.  We say that a string structure on $M$ transgresses to a spin structure on $LM$ (though in this paper we have only discussed isomorphism classes of such structures).  This led to the following conjecture, due independently to H\"ohn and Stolz.

\begin{conj}[H\"ohn--Stolz \cite{Sto96}] Let $M$ be a closed oriented $n$-manifold admitting spin and string structures.  If $M$ admits a metric of positive Ricci curvature, then the Witten genus $\varphi_W(M)=0$.
\end{conj}

Stolz' heuristic argument comes from the hope that there is some Weizenb\"ock-type formula such that positive Ricci curvature on $M$ implies positive scalar curvature on $LM$, which in turn implies $Ker(\Dirac_{LM}) = \Witten (M) = 0$.  Though this line of thinking is far from rigorous, there are no known counterexamples, and the conjecture holds true for homogeneous spaces and complete intersections.  To the author's knowledge, there are no known examples of simply-connected closed manifolds admitting metrics of positive scalar curvature, but not metrics of positive Ricci curvature.  If the conjecture is true, it would provide examples of such manifolds.

Just as $KO$-theory refines the $\Aroof$-genus, there is a cohomology theory $tmf$, or Topological Modular Forms, with string-orientation refining the Witten genus (see \cite{Hop02}):
\[ \xymatrix{& tmf^{-n}(pt) \ar[d] \\
MString^{-n} (pt)\ar[r]^{\Witten } \ar[ur]^{\sigma} & MF_n } \]
The map $tmf^{-*} (pt)\to MF_*$ is a rational isomorphism, but it is not integrally surjective or injective.  In particular, $tmf^{-*}(pt)$ contains a great deal of torsion.  While defining $tmf$ is a subtle process, informally $tmf$ is the universal elliptic cohomology theory, or the elliptic cohomology theory associated to the universal moduli stack of elliptic curves.  Despite several attempts \cite{BDR04, HK04, Seg88, ST04}, there is still no geometric description of $tmf$.  However, it is believed that $tmf$ should provide a natural home for family index theorems on loop spaces.  

One might hope that all the refined invariants in $tmf$ also vanish for string manifolds admitting positive Ricci curvature metrics, giving an analogy of Hitchin's theorem.  However, there exist a fair number of compact non-abelian Lie groups (thus admitting positive Ricci curvature metrics) which are sent to torsion elements in $tmf^{-*}(pt)$ via their left-invariant framing \cite{Hop02}.  For example, in Section \ref{sec:S3} we investigate the case where $M=S^3.$

Conceptually, this is still compatible with the analogy to classical index theory.  The group $Spin(n)$ is a discrete cover of $SO(n)$, so there are no local differences between the bundles $Spin(M)$ and $SO(M)$ and their connections.  However, $Spin^c(n)\to SO(n)$ is an $S^1$-extension, and one must choose a connection on the $S^1$-bundle $Spin^c(M)\to SO(M)$.  The curvature of this connection appears in the Weizenb\"ock formula for the spin$^c$ Dirac operator.  Since $String(M)\to Spin(M)$ has $K(\Z,2)$-fibers, string structures are more analogous to spin$^c$ structures.  When constructing the $S^1$-extension $\widehat{LSpin(M)}\to LSpin(M)$, one really needs an $S^1$-extension with connection \cite{CP89}.  The form $CS_\phalf (g) - \pi^*H_{\cS,g} \in \Omega^3(Spin(M))$ representing $\cS$ transgresses to the curvature (minus a canonically defined term) of this connection on $LSpin(M)$.  One would reasonably expect any Weizenb\"ock-type formula for $\Dirac_{LM}$ to also involve the form $H_{\cS,g}$.  We now ask the following question in an attempt to formulate a connection between $tmf$ and obstructions for certain types of curvature.

\begin{question}\label{question1}Let $M$ be a closed $n$-dimensional manifold with spin structure such that $\phalf(M)=0 \in H^4(M;\Z)$, and let $\cS$ be a specified string class.  Suppose there exists a metric $g$ such that
\[ Ric(g)>0 \quad \text{and} \quad H_{\cS,g}=0\in \Omega^3(M).\]
Does this imply that
\[ \sigma[M,\cS]=0 \in tmf^{-n}(pt) \quad ?\]
\end{question}

\begin{rem}The condition $H_{\cS,g}=0$ for some string class $\cS$ is equivalent to $\check{\phalf}(g)=0\in \check{H}^4(M)$.  This is a strong condition and is not usually satisfied for generic metrics.  While a great deal of information about the characters $\check{\phalf}(g)$ is known for certain manifolds, the author is not aware of any general results guaranteeing the existence or non-existence of such metrics.
\end{rem}

\begin{rem}The condition $H_{\cS,g}=0$ is conformally invariant; if $H_{\cS,g}=0$, then $H_{\cS, e^fg}=0$ for any conformally related metric $e^f g$.  This follows from the conformal invariance of $\check{\phalf}(g)$ and the fact that $0\in \cH^3(M)$ for all metrics.
\end{rem}

We wish to close this discussion by noting that $\Dirac_M$ and $\Dirac_{LM}$ can both be thought of as partition functions of certain 1 and 2-dimensional supersymmetric nonlinear sigma models \cite{Wit99}.  These sigma models require spin and string structures, respectively.  In the 2-dimensional sigma models, the form $H_{\cS,g}$ is used to trivialize the natural connection on a certain determinant line bundle \cite{Wit99, AS02}.  Sometimes, terms in the action of these sigma models are combined and written as the connection $\nabla^{\cS,g}$ discussed in Section \ref{sec:torsionconnections}.

Stolz and Teichner have shown that $KO^{-n}$ is homotopy equivalent to the space of supersymmetric 1-dimensional Euclidean Field Theories of degree $n$ \cite{ST04}, and the spin orientation is (up to homotopy) given by the previously mentioned sigma model.  The hope is that the analogous statement should hold for 2-dimensional field theories with the string orientation $\sigma$ given by these sigma models.  In this context, Question \ref{question1} is essentially asking: if one does not have to add in the terms $H_{\cS,g}$, does positivity of the Ricci curvature imply that the corresponding sigma model is qualitatively trivial?

%%%%%%%%%%%%%%%%%%%%%%%%%%%%%%%%%%%%%%%%%%%%%%%%%%%%
%%%%%%%%%%%%%%%%%%%%%%%%%%%%%%%%%%%%%%%%%%%%%%%%%%%%
\section{Metric connections with torsion}\label{sec:torsionconnections}

We now give an equivalent reformulation of Question \ref{question1} involving the Ricci curvature of a metric connection with torsion.  Given a string class and metric $(\cS, g)$, we define the torsion tensor $T^{\cS,g}$ by
\[T^{\cS,g} \= g^{-1}H_{\cS, g} \in \Omega^1(M; gl(TM)), \]
where $H_{\cS,g}$ is the canonical 3-form from \eqref{eq:Riemannianconstruction}.  This is simply a case of ``raising the indices" and is equivalent to saying $g( T^{\cS,g}_X Y, Z) = H_{\cS, g}(X,Y,Z),$ or in coordinates $T_{ij}^k = g^{rk} H_{ijr}$.  Then
\[ \nabla^{\cS,g} \= \nabla^g + \frac{1}{2} T^{\cS,g}\]
is a metric connection with torsion $T^{\cS,g}$, where $\nabla^g$ denotes the Levi-Civita connection.

In general, torsion tensor $T$ of a connection is called \textit{totally skew-symmetric} if $gT \in \Omega^3(M)$, i.e. $g(T(\cdot, \cdot), \cdot)$ is skew-symmetric in all three variables.  By construction, $\nabla^{\cS,g}$ is a metric connection with totally skew-symmetric torsion.  We also note that the connection $\nabla^{\cS,g}$ still preserves the geodesics of the Levi-Civita connection.  In general for a fixed metric $g$, the the following subsets of connections on $TM$ are equal:
\[ \left\{ \text{Metric connections} \right\} = \left\{ \parbox{1.2in}{Metric connections with $\nabla^g$- geodesics} \right\} = \left\{ \parbox{1.8in}{Metric connections with totally skew-symmetric torsion} \right\} \]
One can easily prove the above by writing any connection $\nabla$ as $\nabla^g + A$ and plugging into the geodesic equation $\nabla_X X=0$ and metric equation $g(\nabla_X Y,Z)=-g(Y,\nabla_X Z)$.

For a torsion-connection $\nabla^T = \nabla^g + \frac{1}{2}T$, we can still define the curvature tensor
\[ R^T_{X,Y}Z \= \left(\nabla^T_X \nabla^T_Y  -\nabla^T_Y \nabla^T_X -\nabla^T_{[X,Y]}\right)Z,\]
and Ricci tensor
\[ Ric^T(X,Y) \= \sum_i g\left( R^T_{e_i, X}Y,e_i \right), \]
where $\{ e_i\}$ is any orthonormal basis.  We let $Ric^g$ denote the Ricci tensor of the Levi-Civita connection.

\begin{lemma}\label{lem:ricci}Suppose that $\nabla^T = \nabla^g + \frac{1}{2}T$ is a metric connection with totally skew-symmetric torsion satisfying $gT = H \in \Omega^3(M)$.  Then the Ricci tensor satisfies
\[Ric^T(X,Y) = Ric^g(X,Y) - \frac{1}{4} \sum_i g(T_{e_i}X, T_{e_i}Y) - \frac{1}{2} d^*H(X,Y).\]
\begin{proof}
Let $\langle \cdot, \cdot \rangle$ denote $g(\cdot,\cdot)$.  Simply expanding using $\nabla^T = \nabla^g + \frac{1}{2} T$, we get
\begin{align*} \langle R^T_{e_i, X}Y, e_i\rangle = & \> \langle \nabla^T_{e_i}\nabla^T_X Y - \nabla^T_X \nabla^T_{e_i}Y -\nabla^T_{[e_i,X]} Y, e_i \rangle \\
= &\> \langle R^g_{e_i, X} Y, e_i \rangle - \frac{1}{4} \langle T_X T_{e_i} Y, e_i \rangle +\frac{1}{2} \langle \nabla^g_{e_i}T_X Y - T_X \nabla^g_{e_i}Y - T_{\nabla^g_{e_i} X} Y, e_i \rangle \\
& + \frac{1}{2} \langle T_{\nabla^g_{e_i} X - \nabla^g_X e_i - [e_i,X]} Y, e_i \rangle \\
=& \> \langle R^g_{e_i, X} Y, e_i \rangle - \frac{1}{4} \langle T_{e_i}X, T_{e_i} Y \rangle + \frac{1}{2} \langle (\nabla^g_{e_i}T)(X,Y), e_i \rangle
\end{align*}
The last term is easily seen to be a tensor.  Using a normal orthonormal frame $\{e_i\}$ at a point (i.e. $\nabla_{e_i}e_j=0$), one easily calculates that
\[ \sum_i  \langle (\nabla^g_{e_i}T)(e_j,e_k), e_i \rangle = \sum_i \partial_i T^i_{jk} = \sum_i \partial_i H_{ijk} = -d^*H (e_j, e_k).  \]
\end{proof}\end{lemma}

The usual Ricci tensor $Ric^g$ is symmetric, and Lemma \ref{lem:ricci} shows that the skew-symmetric part of $Ric^T$ is $-\frac{1}{2}d^*H$.  Since the canonical form $H_{\cS,g}$ satisfies $d^*H_{\cS,g}=0$, it gives rise to a metric connection $\nabla^{\cS,g}$ with symmetric Ricci tensor.  For an arbitrary metric connection $\nabla^T$, we refer to the Ricci curvature $Ric^T(X)\=Ric^T(X,X)$ as the symmetric component of $Ric^T$, which satisfies
\[  Ric^g(X) - Ric^T(X) = \frac{1}{4} \sum_i \| T_{e_i}X \|^2 \geq 0,\]
with equality for all $X$ precisely when $T=0$.  This gives an alternative description of the Levi-Civita connection.
\begin{cor}For a fixed Riemannian metric $g$, the Levi-Civita connection is the unique metric connection which maximizes the Ricci curvature.
\end{cor}

One convenient property of both the Levi-Civita connection and the usual Ricci tensor is the invariance under a global scaling.  A quick check shows that for $\epsilon >0$,
\[
Ric^{\epsilon g}(X) = \sum_i \epsilon g \left( R^{\epsilon g}_{\epsilon^{-\half}e_i, X}X, \epsilon^{-\half}e_i \right) = \sum_i g\left(R^g_{e_i, X}X,e_i \right) = Ric^g(X). 
\]
The form $H_{\cS,g}$ was constructed using a Hodge isomorphism in an adiabatic limit, giving us the scale invariance $H_{\cS, \epsilon g} = H_{\cS, g}$.  However, we use the metric to change $H_{\cS,g}$ into a torsion tensor.  Therefore,
\[ T^{\cS, \epsilon g} = (\epsilon g)^{-1} H_{\cS, \epsilon g} = \epsilon^{-1}T^{\cS, g}.\]
It is more natural then to consider the 1-parameter family of connections $\nabla^{\cS, \epsilon g}$ than any fixed $\nabla^{\cS,g}$.  In the large-volume limit, as $\epsilon \to \infty$, the connection $\nabla^{\cS, \epsilon g}$ converges to the Levi-Civita connection $\nabla^g$.  In the small volume limit, as $\epsilon \to 0$, the terms $T^{\cS, \epsilon g}$ blow up and $\nabla^{\cS, \epsilon g}$ does not converge to a connection unless $H_{\cS,g}=0$.

\begin{question}\label{question2}Let $(M,g,\cS)$ be an $n$-dimensional Riemannian manifold with string class.  Suppose that the Ricci tensor of the modified connection $\nabla^{\cS, g}$ is strictly positive in the small volume scaling limit; that is
\[ \lim_{\epsilon \to 0} Ric(\nabla^{\cS, \epsilon g} )>0.\]
Does this imply $\sigma[M,\cS]=0 \in tmf^{-n}(pt)$?
\end{question}

\begin{prop}\label{prop:equivalenceofhyp} Question \ref{question1} is equivalent to Question \ref{question2}.
\begin{proof}This follows directly from the description of the Ricci tensor in Lemma \ref{lem:ricci}, which implies
\[ Ric^{\epsilon g, \cS} (X) = Ric^g (X) - \frac{1}{4\epsilon} \sum_i \| T_{e_i} X \|^2 .\]
Consequently, if $H_{\cS,g} \neq 0$, then for some $X$,
\[  Ric^{\epsilon g, \cS}(X) \overset{\epsilon \to 0}\longrightarrow -\infty . \]
The simultaneous conditions $Ric(g)>0$ and $H_{\cS,g}=0$ are equivalent to $Ric(\nabla^{\cS, \epsilon g})>0$ for arbitrarily small $\epsilon$.
\end{proof}
\end{prop}

%%%%%%%%%%%%%%%%%%%%%%%%%%%%%%%%%%%%%%%%%%%%%%%%%%%%
%%%%%%%%%%%%%%%%%%%%%%%%%%%%%%%%%%%%%%%%%%%%%%%%%%%%

\section{Homogeneous metrics on $S^3$}\label{sec:S3}

We now examine the canonical 3-forms obtained when $M=S^3$ with a homogeneous metric, and we compare the results with Question \ref{question1}.  We see that Question \ref{question1} has an affirmative answer in this special situation, but it would not if the conditions were weakened.  In particular, there exists a 1-dimensional family of left-invariant metrics $g$ with non-negative Ricci curvature such that the right-invariant framing $\cR$ produces $H_{\cR,g}=0$ and $\sigma[S^3,\cR]\neq 0 \in tmf^{-3}(pt)$.  The previous sentence is also true with left and right swapped.

\subsection{String structures on $S^3$}
Using the isomorphism $S^3 \iso SU(2)\iso Sp(1)$, the left and right-invariant framings induce two string classes which we denote $\cL$ and $\cR$.  The disc $D^4$ inherits a standard framing from its inclusion $D^4 \subset \R^4$, and this restricts to a framing of the stable tangent bundle for $\partial D^4 = S^3$.  We denote the induced string class by $\partial D^4$ and note that, by construction, the string-bordism class $[S^3, \partial D^4] =0 \in MString^{-3}(pt)$.

The set of string classes is a torsor for $H^3(S^3;\Z) \iso \Z$, an affine copy of $\Z$.  In other words, the difference between any two string classes is naturally an integer.  We now determine where the three previously defined string classes live on this affine line, and we use $\Omega c_2 \in H^3(S^3;\Z)$ as our standard generator.  The left and right framings are related by
\[ S^3 \times Spin(3) \overset{L}\longrightarrow Spin(S^3) \overset{R}\longleftarrow S^3 \times Spin(3),\]
and the composition $R^{-1} \of L$ is the Adjoint representation lifted to $Spin$
\[ S^3 \iso SU(2) \overset{Ad}\longrightarrow Spin(su(2)) \iso Spin(3).\]
The difference $\cL - \cR = \pi^* \left( Ad^* \Omega \phalf \right)$.  The Adjoint representation here is an isomorphism of Lie groups and hence an isomorphism on cohomology.  As mentioned in Remark \ref{rem:lowdim}, there is a factor of 2 and minus sign at work: the class $\Omega \phalf$ is twice a generator of $H^3(S^3;\Z)$, and stably $p_1 = -c_2$.  Hence $\Omega \phalf$ is mapped to $-2\Omega c_2$, or $-2\in \Z \iso H^3(S^3;\Z)$, and we use the shorthand
\[ \cL + 2 = \cR.\]

Similarly, we examine the difference between the left-framing and the bounding string structure, and in doing so reference Remark \ref{rem:lowdim}.  The string structure induced from $D^4$ is a framing of the stable tangent bundle.  The normal bundle $\nu \to S^3$ is trivial, and we have the standard isomorphisms of bundles over $S^3$
\[ Spin(TS^3 \oplus \R) \iso Spin(TS^3 \oplus \nu) \iso Spin(D^4) \iso Spin(4).\]
The difference in framing of the two stable bundles differs by the left-multiplication map $S^3 \to Spin(4)$ given by considering $S^3$ as the unit quaternions.  Under the standard isomorphisms $S^3\iso SU(2)$ and $Spin(4) \iso SU(2)\times SU(2)$, this left-multiplication map is the inclusion into the first factor
\[ SU(2) \overset{Id \times \{1\}}\longrightarrow SU(2)\times SU(2) \iso Spin(4). \]
The induced map on cohomology sends $\Omega \phalf$ to $-\Omega c_2$, or $-1 \in \Z \iso H^3(S^3;\Z)$.  Therefore,
\[ \cL + 1 = \partial D^4, \text{ and } \cL + 2 = \partial D^4 + 1 = \cR.\]

The Adams $e$-invariant gives an isomorphism $\pi^s_3 \overset{\iso}\to \Z/24$ and sends the left and right-framings to the two generators (\cite{AS74}).  Our calculations also verify this explicitly.  On a framed $(4k-1)$-dimensional manifold $M$, the $e$-invariant can be computed as follows.  Choose a spin manifold $W$ such that $\partial W=M$ as spin manifolds; such a manifold exists because $MSpin^{4k-1}(pt)=0$.  Using the framing of $TM$, define the Pontryagin classes $p_i(W,M)$ as relative classes in $H^*(W,M)$.  We then obtain $\Aroof(W,M)$ by evaluating $\Aroof(TW, TM)$ on the fundamental class of $W$, where $\Aroof(TW, TM)$ is the $\Aroof$-polynomial with relative Pontryagin classes.  Then,
\[
e[M] =\begin{cases} \Aroof(W,M) \mod \Z, &k \text{ even,}\\
\frac{1}{2} \Aroof(W,M) \mod \Z, &k \text{ odd.} \end{cases}
\]
The $e$-invariant is well-defined as an element of $\Q/\Z$, since choosing a different $W'$ will give $\Aroof(W',M)-\Aroof(W,M) = \Aroof(W'\cup_M(-W)),$ which is an integer (or even integer) by the Atiyah--Singer Index theorem.

If we include metrics so that $(W,\wt{g})$ is a Riemannian spin manifold with boundary $(M,g)$, then we naturally have the Pontryagin forms $p_i(\wt g)\in \Omega^{4k}(W).$

\begin{prop}If $(M,g,\cS)$ is a Riemannian spin 3-manifold with string class, then
\[ e(M, \cS) =  \frac{-1}{48}\int_{W} p_1(\wt{g}) +\frac{1}{24}\int_{M} H_{\cS,g} \mod \Z. \]
\begin{proof}
\[ e[M,\cS] = \frac{1}{2}\int_W \Aroof(W,M) = \frac{1}{2}\int_W \left(1-\frac{1}{24}p_1(W,M) + \cdots \right) = \frac{-1}{48}\int_W p_1(W,M).\]
We now construct a de Rham representative of $p_1(W,M)$.  If $\partial W=M$, then consider the bordism
\[ W \cup_M \left( [0,1] \times M \right) \]
obtained by gluing $\partial W$ to $\{ 0\} \times M$.  The string class $\cS$ gives a stable trivialization $p$ of $Spin(TM)$ up to homotopy, and we let $\Theta_p$ denote the induced flat connection.  Denoting the Levi-Civita connection on $Spin(TM)$ by $\Theta_g$, we have the connection $\Theta(t)$ on $[0,1]\times M$ where
\[ \Theta(t) = t \Theta_p + (1-t) \Theta_g.\]
Finally, define $\wt{\Theta}$ to be the connection on $Spin\left(W \cup_M ( [0,1] \times M) \right)$ induced by $\Theta_{\wt g}$ and $\Theta(t)$.  The form $p_1(\wt \Theta)$ is a de Rham representative of $p_1(W,M)$, and
\begin{align*}\int_{\wt{W}} p_1(\wt \Theta) &= \int_{W} p_1(\wt{g}) + \int_{M^3}\int_{[0,1]} p_1(\Theta(t)). \\
&= \int_W p_1(\wt g) + \int_M CS_{p_1}(\Theta_p, \Theta_g) \\
&= \int_W p_1(\wt g) - 2\int_M CS_{\phalf}(\Theta_g, \Theta_p)
\end{align*}
where $CS_\lambda(\Theta_g,\Theta_p)$ is the general Chern--Simons transgression between two connections.  Lemmas \ref{lem:transgression} and \ref{lem:trivbundle} together imply
\[ \int_M CS_\phalf(\Theta_g, \Theta_p) = \int_M p^* CS(\Theta_g) = \int_M H_{\cS, g}.\]
Therfore,
\begin{align*}
\frac{-1}{48}\int_{\wt{W}} p_1(\wt \Theta) &= \frac{-1}{48}\int_W p_1(\wt g) +\frac{1}{24} \int_M CS_{\phalf}(\Theta_g, \Theta_p) \\
&= \frac{-1}{48}\int_W p_1(\wt g) +\frac{1}{24}\int_M H_{\cS, g}
\end{align*}
\end{proof}\end{prop}

\begin{cor}When $M=S^3$ and $g$ is the standard round metric,
\[ e(S^3, \cS) = \frac{1}{24}\int_{S^3} H_{\cS, g} \mod \Z.\]
\end{cor}

In the next subsection, we calculate $H_{\cS, g}$ for all left-invariant metrics on $S^3$.  Equation \eqref{eq:roundmetric} and the above corollary imply that $e[S^3, \cL]=\frac{-1}{24},$ $e[S^3, \partial D^4]=0,$ and $e[S^3, \cR]=\frac{1}{24}$.  Below is a pictorial description of the space of string classes on $S^3$ and their corresponding string bordism class under $e:MString^{-3}\overset{\iso}\to \Z/24$.
\[ \xymatrix@R=.05cm{ \ar@{<->}[rrrrrrrr]&\bullet & \bullet \ar@{|->}[ddd]&\bullet & \bullet & \bullet & \bullet\ar@{|->}[ddd] & \bullet & \\
&&& \cL \ar@{|->}[dd] & \partial D^4 \ar@{|->}[dd]& \cR\ar@{|->}[dd] \\ \ \\
& &\frac{-2}{24} &\frac{-1}{24} &0 &\frac{1}{24} &\frac{2}{24} } \]

\begin{lemma}\label{lem:transgression}
If $p:M\to P$ is a global section and $\Theta_p$ the induced flat connection, then \[CS_\lambda (\Theta, \Theta_p)  = p^*CS_\lambda(\Theta) \in \Omega^3(M).\]
\begin{proof}
This lemma is essentially a tautology.  Using the notation of \cite{Fre02}, in general $CS_\lambda(\Theta_1, \Theta_0) \= \int_{[0,1]} \lambda(\Theta_t) \in \Omega^{2\bullet -1}(M)$, where $\Theta_t \= t\Theta_1 + (1-t) \Theta_0$ is a connection on $[0,1]\times P \to [0,1]\times M$.  Then, $CS_\lambda(\Theta) \= CS_\lambda(\pi^* \Theta, \Theta_{taut}) \in \Omega^{2\bullet - 1}(P)$, where $\Theta_{taut}$ is the trivial connection induced by the canonical section of $\pi^*P$.  Since one can compute these transgression forms via local frames, and by definition $p^*\Theta_p=0$, we easily see
\begin{align*}
CS_\lambda(\Theta, \Theta_p) &= \int_{[0,1]} \lambda \big( p^* (t\Theta + (1-t)\Theta_p)  \big) = \int_{[0,1]} \lambda \big( p^* t \Theta \big) \\
&= p^* \int_{[0,1]} \lambda (t \Theta) = p^* CS_\lambda(\Theta)
\end{align*}
\end{proof}\end{lemma}

%%%%%%%%%%%%%%%%%%%%%%%%%%%%%%%%%%%%%%%%%%

\subsection{Calculation of canonical 3-forms}
We now investigate Question \ref{question1} by considering left-invariant metrics on $S^3\iso SU(2)$; i.e. metrics $g$ on $SU(2)$ such that left-multiplication is an isometry.  As noted in Proposition \ref{prop:questionS3}, the calculations for right-invariant metrics only differ from those for left-invariant metrics by a sign.  Any such left-invariant metric is determined by its behavior on the tangent space at the identity, so we are considering metrics on the Lie algebra $su(2)$ of left-invariant vector fields.  A global rescaling of $g$ leaves the Ricci tensor and canonical form $H_{\cL, g}$ invariant, hence it does not affect the outcome of Question \ref{question1}.  The space of left-invariant metrics, up to change of oriented basis and global rescaling, is the 2-dimensional space $Sym_{>0}^2(\R^3)/(SO(\R^3) \times \R_+)$, where $Sym^2_{>0}(\R^3)$ denotes the 6-dimensional space of positive-definite $3\times 3$-matrices.

We now give a more computationally explicit description of this space.  Let $\{e_1, e_2, e_3\}$ be the standard basis for $su(2)$ satisfying
\[ [e_1, e_2] = 2 e_3, \> [e_2, e_3] = 2 e_1, \> [e_3, e_1] = 2 e_2.\]
When $\{e_1, e_2, e_3\}$ is an orthonormal basis, the metric is bi-invariant and equal to the standard round metric on $S^3\subset D^4$.  For any $\alpha_1, \alpha_2 \in \R_{>0}$, define the left-invariant metric $g_{\alpha_1, \alpha_2}$ by declaring $\{ \alpha_1 e_1, \alpha_2 e_2, e_3\}$ to be an orthonormal basis.  In the case where $\alpha_2=1$, we recover the 1-parameter family of Berger metrics on $S^3$.  Based on knowledge from \cite{Mil76}, it suffices to consider the 2-parameter family of metrics $\{g_{\alpha_1,\alpha_2}\}$.

\begin{lemma}If $g$ is a left-invariant metric on $SU(2)$, then there exists $\alpha_1, \alpha_2 \in \R_{>0}$ such that $g_{\alpha_1,\alpha_2}$ is isometric to a constant multiple of $g$.
\begin{proof}
Lemma 4.1 in Milnor's \cite{Mil76} implies that for there exists an orthonormal basis $\{E_1, E_2, E_3\}$ for $g$ such that
\[ [E_1, E_2] = \lambda_3 E_3, \quad [E_2, E_3]=\lambda_1 E_1, \quad [E_3, E_1]=\lambda_2 E_2,\]
where $\lambda_i \in \R_{>0}$ (Milnor's $e_i$'s correspond to our $E_i$'s).  For any $(\lambda_1, \lambda_2, \lambda_3)$, it is clear that the orthonormal basis
\[ \{ \frac{\sqrt{\lambda_2\lambda_3}}{2} e_1, \frac{\sqrt{\lambda_3 \lambda_1}}{2} e_2, \frac{\sqrt{\lambda_1 \lambda_2}}{2} e_3\} \]
defines a left-invariant metric isometric to the original $g$.  Finally, we normalize so that the coefficient of $e_3$ is 1.  Hence, there is a surjective map
\begin{align*}
\R^2_{>0} &\longrightarrow \{ \text{Left-inv metrics} \}/ \{ \text{Isom}\times \text{Scale}\} \\
\alpha_1, \alpha_2 &\longmapsto g_{\alpha_1, \alpha_2} 
\end{align*}\end{proof}\end{lemma}

\begin{figure}[htp]
\begin{center}
  \subfigure[Region with positive Ricci curvature]{\includegraphics[scale=.65]{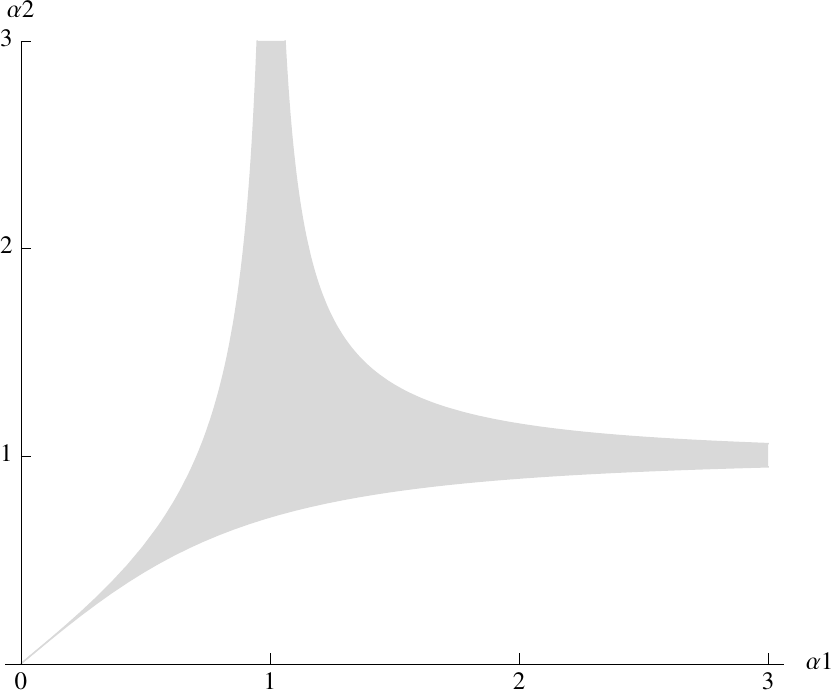}}
  \subfigure[Values of $\int_{S^3}H_{\cL, g_{\alpha_1, \alpha_2}}$]{\includegraphics[scale=.5]{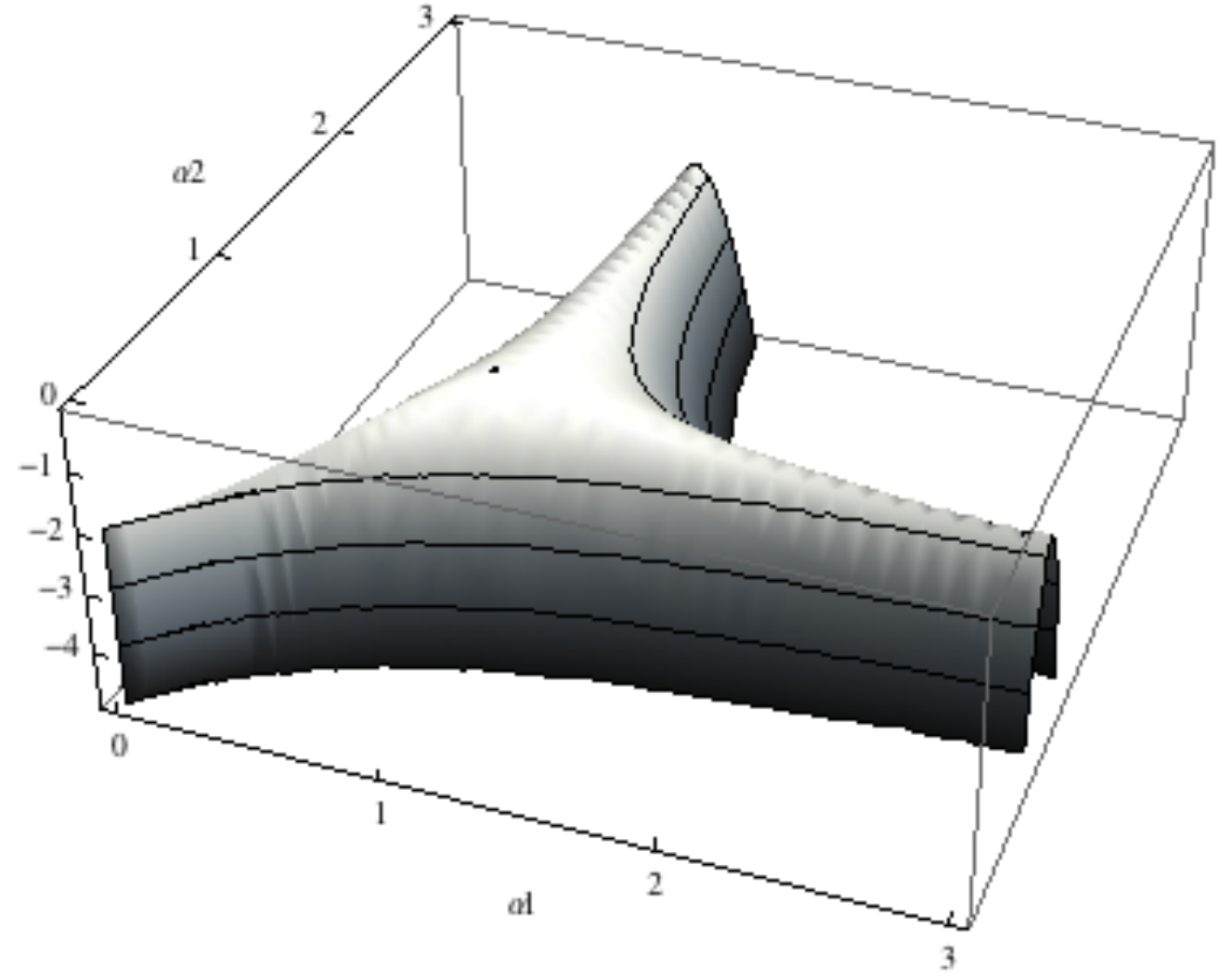}}
\end{center}\caption{}\label{fig1}
\end{figure}

We first calculate the Ricci curvature for $g_{\alpha_1, \alpha_2}$.  A straightforward computation gives the covariant derivative of the Levi-Civita connection in our invariant frame.  The non-zero components are:
\begin{align*}
\langle \nabla_{\alpha_1 e_1} \alpha_2 e_2 , e_3 \rangle &= \alpha_1 \alpha_2 + \frac{\alpha_1}{\alpha_2} - \frac{\alpha_2}{\alpha_1}  \\
\langle \nabla_{\alpha_2 e_2} \alpha_3 e_3, \alpha_1 e_1 \rangle &= \alpha_1 \alpha_2 - \frac{\alpha_1}{\alpha_2} + \frac{\alpha_2}{\alpha_1}  \\
\langle \nabla_{\alpha_3 e_3} \alpha_1 e_1, \alpha_1 e_2 \rangle &= - \alpha_1 \alpha_2 + \frac{\alpha_1}{\alpha_2} + \frac{\alpha_2}{\alpha_1} 
\end{align*}
The Ricci tensor is then diagonalized with eigenvalues
\begin{align*}
Ric(\alpha_1 e_1)&= 2\left(\alpha_1 \alpha_2 - \frac{\alpha_1}{\alpha_2} + \frac{\alpha_2}{\alpha_1}\right) \left(- \alpha_1 \alpha_2 + \frac{\alpha_1}{\alpha_2} + \frac{\alpha_2}{\alpha_1}\right) \\
Ric(\alpha_1 e_2)&= 2\left(- \alpha_1 \alpha_2 + \frac{\alpha_1}{\alpha_2} + \frac{\alpha_2}{\alpha_1}\right) \left(\alpha_1 \alpha_2 + \frac{\alpha_1}{\alpha_2} - \frac{\alpha_2}{\alpha_1} \right) \\
Ric(e_3)&= 2\left(\alpha_1 \alpha_2 + \frac{\alpha_1}{\alpha_2} - \frac{\alpha_2}{\alpha_1} \right) \left(\alpha_1 \alpha_2 - \frac{\alpha_1}{\alpha_2} + \frac{\alpha_2}{\alpha_1}\right).
\end{align*}
Solving inequalities tells us that the Ricci curvature is strictly positive if and only if $(\alpha_1,\alpha_2)$ is in the interior of the region bounded by the three curves
\begin{equation}\label{eq:Riccipos} \alpha_2 = \sqrt{ \frac{\alpha_1^2}{1+\alpha_1^2} }, \quad \alpha_2 = \sqrt{ \frac{\alpha_1^2}{-1+\alpha_1^2} }, \quad \alpha_2 = \sqrt{ \frac{-\alpha_1^2}{-1+\alpha_1^2} }.\end{equation}
This region is shown in Figure \ref{fig1}a.  The Ricci curvature is non-negative with one 0 eigenvalue on the three boundary curves.

Now we calculate the canonical 3-form $H_{\cL, g_{\alpha_1,\alpha_2}} \in \Omega^3(S^3)$.  For dimensional reasons, $H_{\cL, g_{\alpha_1,\alpha_2}}$ is harmonic and therefore
\[ H_{\cL, g_{\alpha_1,\alpha_2}} \in \cH^3(S^3) \iso H^3(S^3;\R) \iso \R\]
with its value in $\R$ determined by integrating over $S^3$.  Lemma \ref{lem:trivbundle} states that we can calculate $H_{\cL, g_{\alpha_1,\alpha_2}}$ by simply calculating the Chern--Simons 3-form $CS_{\phalf}(g_{\alpha_1,\alpha_2})$ on the global frame $\{ \alpha_1 e_1, \alpha_2 e_2, e_3\}$.  This is a straightforward, though lengthy, calculation.

\begin{figure}[htp]
\begin{center}
  \subfigure[Left-invariant metrics]{\includegraphics[scale=.70]{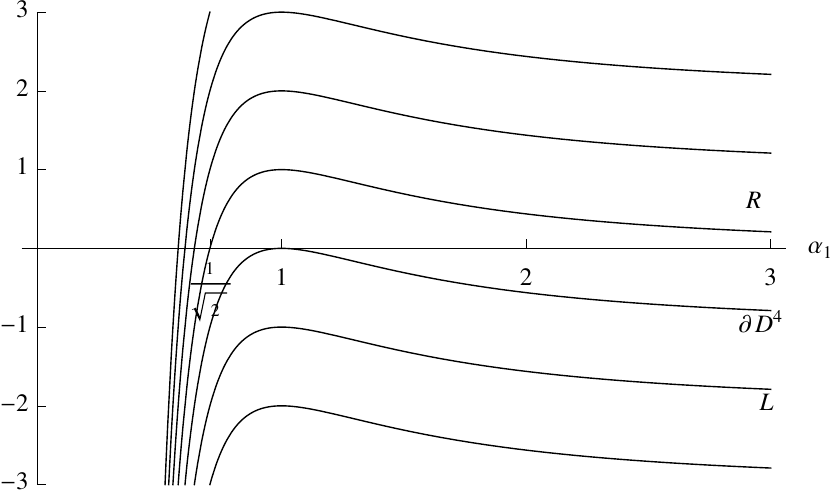}}
  \subfigure[Right-invariant metrics]{\includegraphics[scale=.70]{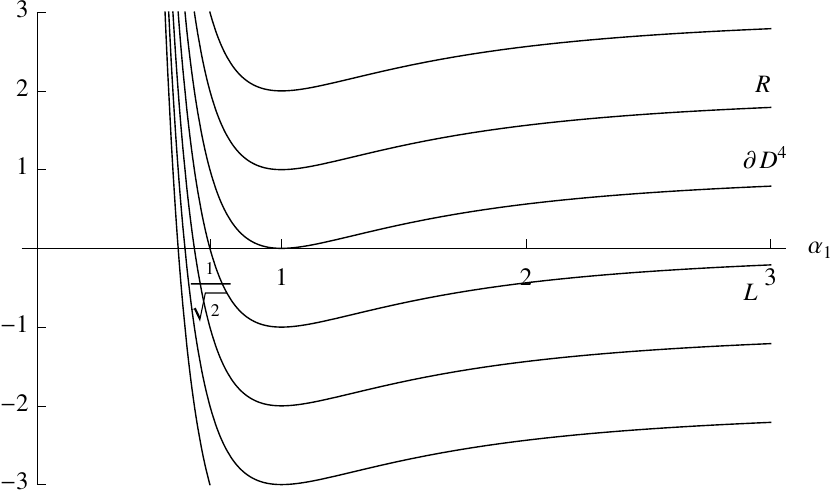}}
\end{center}\caption{$\int_{S^3} H_{\cS, g_{\alpha_1,\alpha_2}}$ on Berger metrics}\label{fig2}
\end{figure}

For the class $\phalf$, the Chern--Simons form is
\[ CS_\phalf(\Theta) = \frac{-1}{16\pi^2} \Tr \left( \Omega\wedge\Theta - \frac{1}{6} \Theta \wedge[ \Theta \wedge \Theta] \right),\]
with $\Tr$ being the ordinary matrix trace.  The normalization constant can be seen from $\phalf(\Theta) =  \frac{-1}{2}c_2(\Theta) = \frac{-1}{2} \frac{1}{8\pi^2} \Tr(\Omega \wedge \Omega)$.  The frame $\{e_i\}$ gives rise to the dual frame $\{ e^i\}$ on $su(2)^*$.  In our global frame, the Chern--Simons form is a constant multiple of $e^1 \wedge e^2 \wedge e^3$, the standard volume form for $SU(2) \iso S^3 \subset D^4$.  Using a direct calculation along with $\int_{S^3}e^1 \wedge e^2 \wedge e^3 = 2\pi^2$, we obtain
\begin{align}\label{eq:HS3} \int_{S^3} H_{\cL, g_{\alpha_1,\alpha_2}} = -\frac{1}{16\pi^2} \int_{S^3}\Tr \left( \Theta \wedge \Omega - \frac{1}{6} \Theta \wedge [\Theta \wedge \Theta] \right)  \\
= -\frac{  \alpha_1^6 \alpha_2^6 - \alpha_1^6 \alpha_2^4  - \alpha_1^4 \alpha_2^6 - \alpha_1^6 \alpha_2^2 - \alpha_1^2 \alpha_2^6
 - \alpha_1^4 \alpha_2^2  -
\alpha_1^2 \alpha_2^4 +  4 \alpha_1^4 \alpha_2^4  +\alpha_1^6 + \alpha_2^6 }{\alpha_1^4 \alpha_2^4}. \nonumber \end{align}
See Figure \ref{fig1}b for a graph of this function.  If we set $\alpha_2=1$ and only consider the usual Berger metrics, we obtain
\begin{equation}\label{eq:berger} \int_{S^3}H_{\cL, g_{\alpha_1, 1}} = -2 + \frac{2\alpha_1^2 - 1}{\alpha_1^4}.\end{equation}
These values are graphed in Figure \ref{fig2}a.  Note that when reduced mod $\Z$, \eqref{eq:berger} coincides with the calculation performed in the original Chern--Simons paper \cite{CS74}.  If we set $\alpha_1=\alpha_2=1$, we obtain the standard bi-invariant metric and see that
\begin{equation}\label{eq:roundmetric}
\int_{S^3}H_{\cL, g_{1,1}} = -1, \quad \int_{S^3}H_{\partial D^4, g_{1,1}} = 0, \quad \int_{S^3}H_{\cR, g_{1,1}} = 1.
\end{equation}

We now analyze \eqref{eq:HS3} on the region $Ric \geq 0$.  The only critical point occurs at $\alpha_1=\alpha_2=1$, where $\int_{S^3} H_{\cL, g_{\alpha_1, \alpha_2}}=-1$ is a maximal value. Furthermore, $\int_{S^3} H_{\cL, g_{\alpha_1, \alpha_2}} = -2$ identically on the three curves bounding the region of positive Ricci curvature.  So, we have the following range of values:
\begin{align*} \left\{ \int_{S^3}H_{\cL, g}\quad | \quad Ric(g)>0, \> g \text{ left-invariant}  \right \} = (-2,-1].\end{align*}
Figure \ref{fig1}b demonstrates this with the help of Mathematica; the level curves for $-2$ are precisely the three functions from \eqref{eq:Riccipos}.

Due to the equivariance of the canonical 3-form under change of string class (see Proposition \ref{prop:Hproperties}), our calculation using $\cL$ gives us $H_{\cS, g_{\alpha_1,\alpha_2}}$ for any other string class $\cS$ by
\[ \int_{S^3} H_{\cL+ j, g_{\alpha_1,\alpha_2}} = j + \int_{S^3} H_{\cL, g_{\alpha_1,\alpha_2}}\]
for any $j \in \Z \iso H^3(S^3;\Z)$.  Therefore,
\begin{align}\label{eq:range} \left\{ \int_{S^3}H_{\cL+j, g}\quad | \quad Ric(g)>0, \> g \text{ left-invariant}  \right \} = (-2+j,-1+j].\end{align}
To graphically demonstrate this, Figure \ref{fig2}a shows the canonical 3-forms for various string classes on the 1-parameter family of left-invariant Berger metrics.  

The entire previous discussion was based on left-invariant Riemannian metrics.  What if we had decided to use right-invariant metrics?  Given an inner product $g_e$ on $T_e SU(2)$, we can form a left-invariant metric $g^L$ and a right-invariant metric $g^R$ by left or right multiplying $g_e$.  The canonical 3-forms are related by the following easy Lemma, whose proof is at the end of this section.

\begin{lemma}\label{lem:leftright}$\displaystyle H_{\cL, g^L} = - H_{\cR, g^R}.$
\end{lemma}
\noindent This fact is graphically demonstrated in Figure \ref{fig2}b below.  In the case of the Berger metrics, note that the Ricci curvature is positive for all $\alpha_1 > \frac{1}{\sqrt{2}}$, and the Ricci curvature is non-negative with a 0 eigenvalue at $\alpha_1 = \frac{1}{\sqrt{2}}$.

\begin{prop}\label{prop:questionS3}If the string class and (left or right)-invariant Riemannian metric $(\cS, g)$ on $S^3$ satisfy
\[ Ric(g) >0, \text{ and } H_{\cS, g} = 0,\]
then $\cS = \partial D^4$ and $g$ is the bi-invariant round metric.  Consequently,
\[ \sigma[S^3, \cS] = 0 \in tmf^{-3}(pt) \iso \Z/24. \]
\begin{proof}
If $g$ is a left-invariant metric with positive Ricci curvature and $H_{\cS, g}=0$, then \eqref{eq:range} implies that $\cS = \cL + 1 = \partial D^4$ with $g$ the bi-invariant metric $g_{1,1}$.

If $g$ is a right-invariant metric, Lemma \ref{lem:leftright} and \eqref{eq:range} imply that
\begin{align} \left\{ \int_{S^3}H_{\cR+j, g}\quad | \quad Ric(g)>0, \> g \text{ right-invariant}  \right \} = [1+j, 2+j).\end{align}
If $H_{\cS,g}=0$, then $\cS = \cR - 1 = \partial D^4$ and $g=g_{1,1}$.  Finally, $[S^3, \partial D^4]=0 \in MString^{-3}$, so $\sigma[S^3, \partial D^4]=0\in tmf^{-3}$.
\end{proof}\end{prop}

We conclude that in this special case, Question \ref{question1} has a very non-trivial affirmative answer.  In particular, there are 1-dimensional families of left and right-invariant metrics which are Ricci non-negative and satisfy $H_{\cR,g}=0$ and $H_{\cL, g}=0$ respectively.  Furthermore, as evidenced by Figure \ref{fig2}, one can find Ricci positive metrics with $H_{\cL,g}$ arbitrarily small but non-zero.  Finally, we point out that for any string class $\cS$, the lift of the Chern--Simons invariant
\[ Met(S^3) \overset{\int H_{\cS, g}}\longrightarrow \R \]
is surjective.  The 1-parameter families of left and right-invariant Berger metrics in Figure \ref{fig2} show this.

\begin{proof}[Proof of Lemma \ref{lem:leftright}]In a left or right-invariant frame, the connection is computed purely in terms of the Lie bracket on vector fields.  On a Lie group $G$,  one can define two Lie algebra structures $[\cdot, \cdot]_L$ and $[\cdot, \cdot]_R$ corresponding to the usual Lie bracket on left or right invariant vector fields.  For $X,Y \in T_eG$, these are related by
\[[X, Y]_L = -[X,Y]_R.\]
If $\Theta_L,\Theta_R$ denote the connections in the two frames, we have $\Theta_L = -\Theta_R$ and $\Omega_L=\Omega_R$, so 
\[ \Tr ( \Theta_L\wedge \Omega_L - \frac{1}{6} \Theta_L \wedge [\Theta_L\wedge \Theta_L] ) = - \Tr ( \Theta_R\wedge \Omega_R - \frac{1}{6} \Theta_R \wedge [\Theta_R \wedge \Theta_R] ) . \] \end{proof}

%%%%%%%%%%%%%%%%%%%%%%%%%%%%%%%%%%

\bibliographystyle{alphanum}
\bibliography{mybibliography}

\begin{thebibliography}{BSCS}

\bibitem[AS1]{AS02}
Orlando Alvarez and I.~M. Singer.
\newblock Beyond the elliptic genus.
\newblock {\em Nuclear Phys. B}, 633(3):309--344, 2002.

\bibitem[AS2]{AS74}
M.~F. Atiyah and L.~Smith.
\newblock Compact {L}ie groups and the stable homotopy of spheres.
\newblock {\em Topology}, 13:135--142, 1974.

\bibitem[BDR]{BDR04}
Nils~A. Baas, Bjorn~Ian Dundas, and John Rognes.
\newblock Two-vector bundles and forms of elliptic cohomology.
\newblock In {\em Topology, geometry and quantum field theory}, volume 308 of
  {\em London Math. Soc. Lecture Note Ser.}, pages 18--45. Cambridge Univ.
  Press, Cambridge, 2004.

\bibitem[BSCS]{BSCS07}
John~C. Baez, Danny Stevenson, Alissa~S. Crans, and Urs Schreiber.
\newblock From loop groups to 2-groups.
\newblock {\em Homology, Homotopy Appl.}, 9(2):101--135, 2007.

\bibitem[CP]{CP89}
R.~Coquereaux and K.~Pilch.
\newblock String structures on loop bundles.
\newblock {\em Comm. Math. Phys.}, 120(3):353--378, 1989.

\bibitem[CS1]{CS85}
Jeff Cheeger and James Simons.
\newblock Differential characters and geometric invariants.
\newblock In {\em Geometry and topology (College Park, Md., 1983/84)}, volume
  1167 of {\em Lecture Notes in Math.}, pages 50--80. Springer, Berlin, 1985.

\bibitem[CS2]{CS74}
Shiing~Shen Chern and James Simons.
\newblock Characteristic forms and geometric invariants.
\newblock {\em Ann. of Math. (2)}, 99:48--69, 1974.

\bibitem[Dai]{Dai91}
Xianzhe Dai.
\newblock Adiabatic limits, nonmultiplicativity of signature, and {L}eray
  spectral sequence.
\newblock {\em J. Amer. Math. Soc.}, 4(2):265--321, 1991.

\bibitem[For]{For95}
Robin Forman.
\newblock Spectral sequences and adiabatic limits.
\newblock {\em Comm. Math. Phys.}, 168(1):57--116, 1995.

\bibitem[Fre1]{Fre95}
Daniel~S. Freed.
\newblock Classical {C}hern-{S}imons theory. {I}.
\newblock {\em Adv. Math.}, 113(2):237--303, 1995.

\bibitem[Fre2]{Fre02}
Daniel~S. Freed.
\newblock Classical {C}hern-{S}imons theory. {II}.
\newblock {\em Houston J. Math.}, 28(2):293--310, 2002.

\bibitem[Hen]{Hen08}
Andr{\'e} Henriques.
\newblock Integrating {$L\sb \infty$}-algebras.
\newblock {\em Compos. Math.}, 144(4):1017--1045, 2008.

\bibitem[Hit]{Hit74}
Nigel Hitchin.
\newblock Harmonic spinors.
\newblock {\em Advances in Math.}, 14:1--55, 1974.

\bibitem[HK]{HK04}
P.~Hu and I.~Kriz.
\newblock Conformal field theory and elliptic cohomology.
\newblock {\em Adv. Math.}, 189(2):325--412, 2004.

\bibitem[Hop]{Hop02}
M.~J. Hopkins.
\newblock Algebraic topology and modular forms.
\newblock In {\em Proceedings of the International Congress of Mathematicians,
  Vol. I (Beijing, 2002)}, pages 291--317, Beijing, 2002. Higher Ed. Press.

\bibitem[Lic]{Lic62}
Andr{\'e} Lichnerowicz.
\newblock Laplacien sur une vari\'et\'e riemannienne et spineurs.
\newblock {\em Atti Accad. Naz. Lincei Rend. Cl. Sci. Fis. Mat. Nat. (8)},
  33:187--191, 1962.

\bibitem[Liu]{Liu96}
Kefeng Liu.
\newblock Modular forms and topology.
\newblock In {\em Moonshine, the {M}onster, and related topics ({S}outh
  {H}adley, {MA}, 1994)}, volume 193 of {\em Contemp. Math.}, pages 237--262.
  Amer. Math. Soc., Providence, RI, 1996.

\bibitem[LM]{LM89}
H.~Blaine Lawson, Jr. and Marie-Louise Michelsohn.
\newblock {\em Spin geometry}, volume~38 of {\em Princeton Mathematical
  Series}.
\newblock Princeton University Press, Princeton, NJ, 1989.

\bibitem[Mil]{Mil76}
John Milnor.
\newblock Curvatures of left invariant metrics on {L}ie groups.
\newblock {\em Advances in Math.}, 21(3):293--329, 1976.

\bibitem[MM]{MM90}
Rafe~R. Mazzeo and Richard~B. Melrose.
\newblock The adiabatic limit, {H}odge cohomology and {L}eray's spectral
  sequence for a fibration.
\newblock {\em J. Differential Geom.}, 31(1):185--213, 1990.

\bibitem[Red1]{Redden06}
Corbett Redden.
\newblock {\em Canonical metric connections associated to string structures}.
\newblock PhD thesis, University of Notre Dame, 2006.

\bibitem[Red2]{Redden08p1}
Corbett Redden.
\newblock Harmonic forms on principal bundles, 2008.
\newblock [arXiv.dg:0810.4578].

\bibitem[Seg]{Seg88}
Graeme Segal.
\newblock Elliptic cohomology (after {L}andweber-{S}tong, {O}chanine, {W}itten,
  and others).
\newblock {\em Ast\'erisque}, (161-162):Exp.\ No.\ 695, 4, 187--201 (1989),
  1988.

\bibitem[Ser]{Ser51}
Jean-Pierre Serre.
\newblock Homologie singuli\`ere des espaces fibr\'es. {A}pplications.
\newblock {\em Ann. of Math. (2)}, 54:425--505, 1951.

\bibitem[SP]{SP09}
Christopher Schommer-Pries.
\newblock A finite-dimensional string 2-group, 2009.
\newblock [arXiv:0911.2483].

\bibitem[SS]{SS08a}
James Simons and Dennis Sullivan.
\newblock Axiomatic characterization of ordinary differential cohomology.
\newblock {\em J. Topol.}, 1(1):45--56, 2008.

\bibitem[ST]{ST04}
Stephan Stolz and Peter Teichner.
\newblock What is an elliptic object?
\newblock In {\em Topology, geometry and quantum field theory}, volume 308 of
  {\em London Math. Soc. Lecture Note Ser.}, pages 247--343. Cambridge Univ.
  Press, Cambridge, 2004.

\bibitem[Sto1]{Sto92}
Stephan Stolz.
\newblock Simply connected manifolds of positive scalar curvature.
\newblock {\em Ann. of Math. (2)}, 136(3):511--540, 1992.

\bibitem[Sto2]{Sto96}
Stephan Stolz.
\newblock A conjecture concerning positive {R}icci curvature and the {W}itten
  genus.
\newblock {\em Math. Ann.}, 304(4):785--800, 1996.

\bibitem[Wal]{Wal09}
Konrad Waldorf.
\newblock String connections and chern-simons theory, 2009.
\newblock [arXiv:0906.0117].

\bibitem[Wit1]{Wit88}
Edward Witten.
\newblock The index of the {D}irac operator in loop space.
\newblock In {\em Elliptic curves and modular forms in algebraic topology
  (Princeton, NJ, 1986)}, volume 1326 of {\em Lecture Notes in Math.}, pages
  161--181. Springer, Berlin, 1988.

\bibitem[Wit2]{Wit99}
Edward Witten.
\newblock Index of {D}irac operators.
\newblock In {\em Quantum fields and strings: a course for mathematicians, Vol.
  1, 2 (Princeton, NJ, 1996/1997)}, pages 475--511. Amer. Math. Soc.,
  Providence, RI, 1999.

\end{thebibliography}

\end{document}